\documentclass[11pt]{amsart}
\usepackage[T1]{fontenc}
\usepackage{lmodern}
\usepackage{amsmath,amssymb,amsthm,mathtools}
\usepackage[margin=1in]{geometry}
\usepackage{microtype}
\usepackage[hidelinks]{hyperref}
\hypersetup{pdftitle={Deformations and Local Geometry of Lie Triple System Morphisms}}
\allowdisplaybreaks[2]
\numberwithin{equation}{section}
\newtheorem{theorem}{Theorem}[section]
\newtheorem{proposition}[theorem]{Proposition}
\newtheorem{lemma}[theorem]{Lemma}
\newtheorem{corollary}[theorem]{Corollary}
\theoremstyle{definition}
\newtheorem{definition}[theorem]{Definition}
\newtheorem{example}[theorem]{Example}
\newtheorem{remark}[theorem]{Remark}
\newcommand{\K}{\mathbb K}
\newcommand{\R}{\mathbb R}
\newcommand{\g}{\mathfrak g}
\newcommand{\h}{\mathfrak h}
\newcommand{\s}{\mathfrak s}
\newcommand{\m}{\mathfrak m}
\newcommand{\C}{\mathcal C}
\newcommand{\Hom}{\operatorname{Hom}}
\newcommand{\End}{\operatorname{End}}

\newcommand{\InnDer}{\operatorname{InnDer}}
\newcommand{\Inn}{\operatorname{Inn}}
\newcommand{\im}{\operatorname{im}}
\newcommand{\id}{\operatorname{id}}
\newcommand{\pr}{\operatorname{pr}}
\newcommand{\Sh}{\operatorname{Sh}}
\newcommand{\Gr}{\operatorname{Gr}}
\newcommand{\tr}{\operatorname{tr}}
\newcommand{\ad}{\operatorname{ad}}

\newcommand{\B}{\mathsf B}
\title[Deformations and local geometry of morphisms]{Deformations and Local Geometry\\ of Lie Triple System Morphisms}
\author{Wen Teng}
\address{School of Mathematics and Statistics, Guizhou University of Finance
and Economics, Guiyang 550025, P. R. of China}
\email{tengwen@mail.gufe.edu.cn}

\author{Ripan Saha}
\address{Department of Mathematics, Raiganj University
Raiganj 733134, West Bengal, India}
\email{ripanjumaths@gmail.com}
\date{}

\begin{document}
\begin{abstract}
In this paper, first we study the deformations, stability and rigidity of Lie triple system morphisms by means of higher derived brackets and augmented Yamaguti cohomology. The vanishing of the first cohomology gives formal and local inner rigidity. The vanishing of the second cohomology removes the formal extension obstructions and, in finite dimension over $\R$, gives local smoothness and stability under simultaneous perturbations of the two products. Next, we establish the relationship with Lie algebra homomorphisms in each part of the theory. The universal enveloping Lie algebra identifies a Lie triple morphism with a grading-preserving Lie homomorphism. We obtain the corresponding first cohomology isomorphism, the equivalence of formal and smooth deformations, and a stability result for graded Lie bracket perturbations. For Lie triple products induced by binary Lie brackets, we construct a comparison of differentials and obstruction classes. Under perfectness and a zero-centralizer condition, the local morphism spaces coincide. Finally, we relate triple subsystem graphs to graded Lie subalgebra graphs and compare their normal differentials.
\end{abstract}
\maketitle
\enlargethispage{2pt}

\section{Introduction}

Lie triple systems describe the infinitesimal structure of symmetric spaces. Their algebraic structure was developed by Jacobson and Lister \cite{Jacobson,Lister}. In the finite-dimensional real setting, the correspondence between Lie triple systems and germs of pointed symmetric spaces identifies a triple morphism with the tangent map of a local symmetric-space morphism \cite{Smirnov}. Thus a deformation of a triple morphism describes a variation of the infinitesimal data of a map preserving the local symmetric-space structure. This raises geometric questions about whether nearby maps are equivalent under inner symmetries, whether infinitesimal variations extend to actual families, and whether a map persists when the structures on its source and target change.

For fixed triple products $\pi$ and $\mu$ on finite-dimensional real vector spaces $\g$ and $\h$, these questions concern the local structure of a polynomial zero set in $\Hom(\g,\h)$. The morphism space $\operatorname{Mor}(\g,\h)$ is the zero set of the defect
\[
 S_{\pi,\mu}(q)(x,y,z)=\mu(qx,qy,qz)-q\pi(x,y,z).
\]
Its local geometry at a morphism $f$ involves the smoothness of this zero set, the position of the inner automorphism orbit through $f$, and the persistence of nearby zeros as $\pi$ and $\mu$ vary. This provides a finite-dimensional setting in which the relation between the morphism equation and its linearization can be studied directly.

Cohomology gives the infinitesimal description of this geometry. The derivative of the defect at $f$ is the first differential $\delta_f^1$, so $Z_f^1=\ker\delta_f^1$ is the space of solutions of the linearized morphism equation. The degree-zero differential records infinitesimal inner actions, and its image $B_f^1$ is the tangent space to the inner orbit. Consequently, the augmented group $H_f^1=Z_f^1/B_f^1$ describes first-order variations modulo inner equivalence. Higher coefficients of a formal deformation give obstruction classes in $H_f^2$. We study smooth families of morphisms through the local geometry of the morphism equation.

Yamaguti introduced a cohomology theory intrinsic to Lie triple systems \cite{Yamaguti}. Harris studied a related theory using Lie algebras with involution \cite{Harris}. The deformation theory of triple products was developed by Kubo and Taniguchi \cite{KuboTaniguchi}, and Zhang explained Yamaguti cohomology through Leibniz cohomology \cite{Zhang}. Xia, Sheng and Tang constructed a graded Lie algebra that controls Lie triple systems and recovers the Yamaguti differential \cite{XiaShengTang}. We use this graded Lie algebra with its last-three-variable symmetry conditions throughout the paper.

The derived bracket method of Voronov \cite{Voronov}, together with twisting of $L_\infty$-algebras \cite{Getzler}, gives a convenient description of equations for linear maps. Fr\'egier and Zambon applied this method to simultaneous deformations of algebras and morphisms \cite{FregierZambon}. In the setting of Lie triple systems, Hajjaji constructed the controlling algebra and cohomology of crossed homomorphisms \cite{Hajjaji}. For the zero action, crossed homomorphisms are ordinary morphisms. Thus the Maurer--Cartan description used below is a specialization of that construction. Chtioui, Hajjaji, Mabrouk and Makhlouf developed a cohomology and deformation theory for $\mathcal O$-operators on Lie triple systems, using a Lie $3$-algebra whose Maurer--Cartan elements are precisely such operators \cite{ChtiouiHajjajiMabroukMakhlouf}. Cohomology for simultaneous deformations of Lie--Yamaguti algebra morphisms was studied by Mondal and Saha \cite{MondalSaha}. The positive-degree complex considered here keeps the two Lie triple systems fixed; varying their products is treated separately as a parameter-dependent stability problem.

Our geometric approach follows the rigidity and stability framework for Lie algebras described by Crainic, Sch\"atz and Struchiner \cite{CrainicSchaetzStruchiner} and the corresponding study of $3$-Lie algebra morphisms by Jiang, Sheng and Sun \cite{JiangShengSun}. Two features of Lie triple systems enter the argument. First, the triple product is skew in its first two arguments and satisfies a cyclic identity. We use the corresponding cochain spaces and coefficient representation. Second, the degree-zero term records inner equivalence. This term gives the first boundary space and enters the comparison with the normal complex of a graph subsystem.

With these conventions, our main results give cohomological criteria for the local structure of the morphism space. If $H_f^1=0$, every sufficiently close morphism belongs to the inner orbit of $f$, and the required inner automorphism can be chosen smoothly near $f$. Thus the orbit contains a neighborhood of $f$ in the morphism locus. If $H_f^2=0$, that locus is a smooth submanifold near $f$, with
\[
 T_f\operatorname{Mor}(\g,\h)=Z_f^1,
 \qquad \dim_f\operatorname{Mor}(\g,\h)=\dim Z_f^1.
\]
Every first cocycle is then realized as the derivative of a smooth family of morphisms. These assertions identify the local dimension and distinguish orbit directions from the infinitesimal directions measured by $H_f^1$. We prove them by an elementary exactness lemma and an explicit compatibility identity for the defect. The identity accounts for the relations among the morphism equations, and the vanishing of $H_f^2$ supplies the exactness needed for local smoothness.

The stability result describes how this geometry depends on the two triple products. Under the same condition $H_f^2=0$, every sufficiently small simultaneous perturbation of $\pi$ and $\mu$ through Lie triple system products admits a morphism close to $f$. Moreover, this morphism can be chosen as the restriction of a smooth function of the product tensors. Hence any smooth family of such product pairs admits a compatible smooth family of morphisms through $f$. The smooth choice is defined on the ambient tensor space and then restricted to the locus of Lie triple system products. If $f$ is injective, surjective or invertible, the nearby morphisms retain that property after restricting the neighborhood.

Triple subsystems give a second geometric realization of these results. The $k$-dimensional triple subsystems of a finite-dimensional real Lie triple system $E$ form a locus in the Grassmannian $\Gr_k(E)$. At a subsystem $\s$, the ambient tangent space is $\Hom(\s,E/\s)$, and the normal cohomology describes the linearized closure condition and its obstructions. We obtain corresponding criteria for local inner rigidity, smoothness and stability. In particular, the graph $G_f=\{(x,fx):x\in\g\}$ identifies the local morphism space near $f$ with the local subsystem space near $G_f$ in $\Gr_{\dim\g}(\g\oplus\h)$. The comparison of their complexes, including the degree-zero boundaries, shows that deformations of maps and deformations of their graph subsystems have the same local geometry and the same infinitesimal inner equivalence.

A further purpose of this paper is to relate this local geometry to that of Lie algebra homomorphisms. The universal enveloping Lie algebra of a Lie triple system is an ordinary Lie algebra with a $\mathbb Z_2$-grading. Jacobson's universal embedding and Smirnov's graded formulation identify Lie triple morphisms with grading-preserving Lie homomorphisms \cite{Jacobson,Smirnov}. We apply this established correspondence to first cohomology and formal deformations. In finite dimension, polynomial extension and linear restriction identify the local morphism spaces and preserve smooth deformations and inner orbits. Local smoothness and rigidity therefore have equivalent interpretations in these two settings. The recognition criterion for universal envelopes also gives stability under small graded Lie bracket perturbations. The Lie equivalences in this comparison are generated by even inner derivations.

There is another comparison when the source and target triple products are induced by prescribed binary Lie brackets through $[x,y,z]=[[x,y],z]$. A Lie homomorphism always preserves the associated triple products, but a Lie triple morphism need not preserve the binary bracket. We express the triple morphism defect in terms of the binary defect, thereby relating the first differentials and the second obstruction classes. If the source Lie algebra is perfect and the centralizer of its image in the target is zero, the comparison operator is injective. The formal morphism equations and the local morphism loci then agree near the given homomorphism. Perfectness of the target further identifies the inner boundary spaces and the local inner rigidity conditions. These assumptions specify when the local geometry for prescribed binary Lie brackets is recovered from the ternary problem.

The paper is organized as follows. In Sect.~2, we recall the controlling graded Lie algebra and give the Maurer--Cartan characterization, the universal Lie correspondence and the binary defect identity. In Sect.~3, we study augmented cohomology and formal obstructions, and compare them with the Chevalley--Eilenberg theory of Lie homomorphisms. In Sect.~4, we prove local rigidity and stability, together with their graded Lie and binary Lie consequences. In Sect.~5, we study triple subsystems, their normal complexes and their relationship with Lie subalgebra graphs.

Unless otherwise stated, vector spaces are over a field $\K$ of characteristic zero. The algebraic and formal results hold in arbitrary dimension over $\K$. In Sections~4 and~5, all vector spaces are finite-dimensional over $\R$ whenever smoothness, neighborhoods or stability are discussed. All topologies there are the usual finite-dimensional topologies. The algebraic graph comparison also holds in arbitrary dimension.

\section{Derived brackets and the morphism equation}

In this section, we recall the graded Lie algebra for Lie triple systems and apply higher derived brackets to the morphism equation. Subsequently, we establish its relationship with graded Lie homomorphisms through universal envelopes. For triple products induced by binary Lie brackets, we compare the two morphism defects and give an injectivity criterion.

\begin{definition}
A Lie triple system is a vector space $L$ with a trilinear product satisfying
\begin{align}
 [x,x,y]&=0, \label{eq:lts-skew}\\
 [x,y,z]+[y,z,x]+[z,x,y]&=0, \label{eq:lts-cyclic}\\
 [x,y,[u,v,w]]&=[[x,y,u],v,w]+[u,[x,y,v],w]
                    +[u,v,[x,y,w]]. \label{eq:lts-fi}
\end{align}
A morphism $f:(\g,\pi)\longrightarrow(\h,\mu)$ is a linear map such that
\begin{equation}\label{eq:morphism}
 \mu(fx,fy,fz)=f\pi(x,y,z)
 \qquad(x,y,z\in\g).
\end{equation}
\end{definition}

For $x,y\in L$, the map $\ad_{x\wedge y}z=[x,y,z]$ is a derivation by \eqref{eq:lts-fi}. These maps span a Lie algebra $\InnDer(L)$, since
\begin{equation}\label{eq:inner-lie}
 [\ad_{x\wedge y},\ad_{u\wedge v}]
 =\ad_{[x,y,u]\wedge v}+\ad_{u\wedge[x,y,v]}.
\end{equation}
The notation $\ad_X$ extends linearly to $X\in\wedge^2L$ and defines a surjective linear map $\wedge^2L\to\InnDer(L)$.

For vector spaces $L,V$, set $\C^0(L,V)=\Hom(L,V)$. For $p\geq1$, let $\C^p(L,V)$ consist of the maps
\[
 P:(L\otimes L)^{\otimes p}\otimes L\longrightarrow V
\]
that satisfy
\begin{gather}
 P(X_1,\ldots,X_{p-1},x,x,y)=0,\label{eq:cochain-skew} \\
 P(X_1,\ldots,X_{p-1},x,y,z) +P(X_1,\ldots,X_{p-1},y,z,x) +P(X_1,\ldots,X_{p-1},z,x,y)=0.\label{eq:cochain-cyclic}
\end{gather}
Here $X_i\in L\otimes L$, and the last three entries are written separately. These are the cochain conventions of \cite{ChtiouiHajjajiMabroukMakhlouf,XiaShengTang}. The cochain symmetries are exactly \eqref{eq:cochain-skew}--\eqref{eq:cochain-cyclic}, imposed on the last three entries; the preceding arguments remain tensor pairs.

We use the following insertion convention. If $P\in\C^p(L,L)$, $Q\in\C^q(L,L)$ and $X_i=x_i\otimes y_i$, define, for $1\leq k\leq p$,
\begin{align}
 &(P\circ_k Q)(X_1,\ldots,X_{p+q},z) =\sum_{\sigma\in\Sh(k-1,q)}\operatorname{sgn}(\sigma) P\bigl(X_{\sigma(1)},\ldots,X_{\sigma(k-1)}, \notag \\
 &\qquad Q(X_{\sigma(k)},\ldots,X_{\sigma(k+q-1)},x_{k+q})\otimes y_{k+q}, X_{k+q+1},\ldots,X_{p+q},z\bigr) \notag \\
 &\quad+\sum_{\sigma\in\Sh(k-1,q)}\operatorname{sgn}(\sigma) P\bigl(X_{\sigma(1)},\ldots,X_{\sigma(k-1)}, \notag \\
 &\qquad x_{k+q}\otimes Q(X_{\sigma(k)},\ldots,X_{\sigma(k+q-1)},y_{k+q}), X_{k+q+1},\ldots,X_{p+q},z\bigr).\label{eq:insert-pair}
\end{align}
In this formula an unshuffle permutes the first $k+q-1$ indices. The remaining insertion is
\begin{align}
 &(P\circ_{p+1}Q)(X_1,\ldots,X_{p+q},z) \notag\\
 &=\sum_{\sigma\in\Sh(p,q)}\operatorname{sgn}(\sigma)
 P\bigl(X_{\sigma(1)},\ldots,X_{\sigma(p)},
 Q(X_{\sigma(p+1)},\ldots,X_{\sigma(p+q)},z)\bigr). \label{eq:insert-last}
\end{align}
Empty strings of arguments are omitted. Put
\begin{equation}\label{eq:graded-bracket}
 P\circ Q=\sum_{k=1}^{p+1}(-1)^{(k-1)q}P\circ_kQ,
 \qquad [P,Q]_{\B}=P\circ Q-(-1)^{pq}Q\circ P.
\end{equation}

\begin{proposition}[\cite{ChtiouiHajjajiMabroukMakhlouf,XiaShengTang}]\label{prop:control}
The bracket \eqref{eq:graded-bracket} makes $\C^\bullet(L,L)$ a graded Lie algebra. A map $\pi\in\C^1(L,L)$ is a Lie triple system product if and only if $[\pi,\pi]_{\B}=0$.
\end{proposition}

For a degree-one cochain satisfying \eqref{eq:cochain-skew}--\eqref{eq:cochain-cyclic}, the Maurer--Cartan equation is precisely \eqref{eq:lts-fi}.

We use the symmetric convention for $L_\infty[1]$-algebras: the maps $\ell_k$ have degree $1$, are graded symmetric, and satisfy
\[
 \sum_{i=1}^{n}\ \sum_{\sigma\in\Sh(i,n-i)}
 \epsilon(\sigma)\,
 \ell_{n-i+1}\bigl(\ell_i(a_{\sigma(1)},\ldots,a_{\sigma(i)}),
 a_{\sigma(i+1)},\ldots,a_{\sigma(n)}\bigr)=0.
\]
Here $\epsilon(\sigma)$ is the Koszul sign. A degree-zero element $a$ is Maurer--Cartan when $\sum_{k\geq1}\ell_k(a,\ldots,a)/k!=0$. All such sums in the present morphism problem are finite.

Let $G=\g\oplus\h$ and $\mathcal L=\C^\bullet(G,G)$. A map $a\in\C^p(\g,\h)$ is extended to $G$ by taking the $\g$-component of every input and putting its value in the $\h$-component. Denote this extension by $\widehat a$. It vanishes as soon as an input lies in $\h$. Identify
\[
 \mathcal A=\bigoplus_{p\geq0}\C^p(\g,\h)
\]
with the resulting subspace of $\mathcal L$. Let $\mathsf P:\mathcal L\to\mathcal A$ be restriction to $\g$ inputs followed by projection to $\h$. The direct product law is
\[
 \Delta\bigl((x,u),(y,v),(z,w)\bigr)
       =\bigl(\pi(x,y,z),\mu(u,v,w)\bigr).
\]

\begin{proposition}\label{prop:derived}
The quadruple $(\mathcal L,\mathcal A,\mathsf P,\Delta)$ is a V-data. Its higher derived brackets are
\begin{align}
 \ell_1(a)&=\mathsf P[\pi,\widehat a]_{\B}, \label{eq:l1}\\
 \ell_3(a,b,c)&=\mathsf P[[[\mu,\widehat a]_{\B},\widehat b]_{\B},\widehat c]_{\B}, \label{eq:l3}\\
 \ell_2&=0,\qquad \ell_k=0\quad(k\geq4). \label{eq:vanish}
\end{align}
They define an $L_\infty[1]$-algebra on $\mathcal A$.
\end{proposition}

\begin{proof}
The insertion of an element of $\mathcal A$ into another element of $\mathcal A$ is zero, since the inserted value belongs to $\h$ and the outer cochain vanishes whenever one of its inputs belongs to $\h$. Thus $\mathcal A$ is abelian. A cochain belongs to $\ker\mathsf P$ precisely when its value on inputs from $\g$ belongs to $\g$. This property is preserved by insertions, and hence $\ker\mathsf P$ is a graded Lie subalgebra. The two triple products act on different summands, so
\[
 [\pi+\mu,\pi+\mu]_{\B}
 =[\pi,\pi]_{\B}+2[\pi,\mu]_{\B}+[\mu,\mu]_{\B}=0,
 \qquad \mathsf P(\pi+\mu)=0.
\]
Voronov's theorem \cite{Voronov} therefore gives
\[
 \ell_k(a_1,\ldots,a_k)
 =\mathsf P[\cdots[[\pi+\mu,\widehat a_1]_{\B},
                  \widehat a_2]_{\B},\ldots,\widehat a_k]_{\B}.
\]
For the source product, $[\pi,\widehat a_1]_{\B}$ has inputs in $\g$ and output in $\h$, so its bracket with a second extended cochain vanishes. For the target product, all three inputs of $\mu$ must be supplied by extended cochains before projection to $\mathcal A$ can be nonzero. Thus only the third derived bracket receives a contribution from $\mu$, and inserting a fourth extended cochain gives zero. This proves \eqref{eq:l1}--\eqref{eq:vanish}.
\end{proof}

\begin{theorem}\label{thm:MC}
A linear map $f:\g\to\h$ is a Lie triple system morphism if and only if it is Maurer--Cartan in $(\mathcal A,\ell_1,\ell_3)$.
\end{theorem}

\begin{proof}
For $x,y,z\in\g$, the insertion formulas give
\[
 \ell_1(f)(x,y,z)=-f\pi(x,y,z),\qquad
 \ell_3(f,f,f)(x,y,z)=6\mu(fx,fy,fz).
\]
Consequently, the Maurer--Cartan equation evaluated at $x,y,z$ is
\[
 -f\pi(x,y,z)+\frac{1}{6}\,6\mu(fx,fy,fz)=0.
\]
This is precisely the morphism identity \eqref{eq:morphism}.
\end{proof}

If $f$ is a morphism, the twisted brackets \cite{Getzler} are
\begin{equation}\label{eq:twist}
 \ell_1^f(a)=\ell_1(a)+\tfrac12\ell_3(f,f,a),\quad
 \ell_2^f(a,b)=\ell_3(f,a,b),\quad
 \ell_3^f=\ell_3.
\end{equation}
Consequently, $f+a$ is a morphism if and only if
\begin{equation}\label{eq:twisted-MC}
 \ell_1^f(a)+\tfrac12\ell_2^f(a,a)+\tfrac16\ell_3^f(a,a,a)=0.
\end{equation}
For linear maps $a,b$, the quadratic bracket has the explicit form
\begin{align*}
 \ell_2^f(a,b)(x,y,z)
 &=\mu(ax,by,fz)+\mu(bx,ay,fz)
   +\mu(ax,fy,bz)+\mu(bx,fy,az)\\
 &\quad+\mu(fx,ay,bz)+\mu(fx,by,az).
\end{align*}
For three linear maps $a_1,a_2,a_3$, the cubic bracket is
\[
 \ell_3(a_1,a_2,a_3)(x,y,z)
 =\sum_{\sigma\in S_3}\mu(a_{\sigma(1)}x,a_{\sigma(2)}y,a_{\sigma(3)}z).
\]

For any linear map $q$, define its defect and a linear operator by
\begin{align}
 S_{\pi,\mu}(q)(x,y,z)&=\mu(qx,qy,qz)-q\pi(x,y,z), \label{eq:defect}\\
 \mathcal B_{\pi,\mu,q}(b)&=\ell_1(b)+\tfrac12\ell_3(q,q,b),
 \qquad b\in\C^1(\g,\h). \label{eq:compatibility-map}
\end{align}
The value of $\mathcal B_{\pi,\mu,q}$ lies in $\C^2(\g,\h)$.

\begin{lemma}\label{lem:bianchi}
For every $q\in\Hom(\g,\h)$,
\begin{equation}\label{eq:bianchi}
 \mathcal B_{\pi,\mu,q}\bigl(S_{\pi,\mu}(q)\bigr)=0.
\end{equation}
\end{lemma}

\begin{proof}
The extension of $q$ vanishes on $\h$, so $\widehat q^{\,2}=0$ and $(\id+\widehat q)^{-1}=\id-\widehat q$. Conjugating the direct product law by $\id+\widehat q$ gives the trilinear map
\[
 (\id-\widehat q)\Delta
       (\id+\widehat q,\id+\widehat q,\id+\widehat q).
\]
Its value on $(x,u),(y,v),(z,w)$ is
\[
 \bigl(\pi(x,y,z),
       \mu(u+qx,v+qy,w+qz)-q\pi(x,y,z)\bigr).
\]
Conjugation preserves the Lie triple identities. Since $\ker\mathsf P$ is a graded Lie subalgebra and $\mathcal A$ is abelian, projecting the graded self-bracket of this conjugated law gives
\begin{equation*}
 0={}2\mathsf P\Bigl[ (\id-\widehat q)\Delta (\id+\widehat q,\id+\widehat q,\id+\widehat q),\qquad \widehat{\mu(q\,\cdot,q\,\cdot,q\,\cdot) -q\pi(\cdot,\cdot,\cdot)}\Bigr]_{\B}.
\end{equation*}
For every $b\in\C^1(\g,\h)$, the insertion formulas and \eqref{eq:l1}--\eqref{eq:vanish} give
\begin{equation*}
 \mathsf P\Bigl[ (\id-\widehat q)\Delta (\id+\widehat q,\id+\widehat q,\id+\widehat q), \widehat b\Bigr]_{\B} =\mathsf P[\pi,\widehat b]_{\B} +\tfrac12\mathsf P[[[\mu,\widehat q]_{\B},\widehat q]_{\B}, \widehat b]_{\B}.
\end{equation*}
Substituting the full trilinear expression
$\mu(q\,\cdot,q\,\cdot,q\,\cdot)-q\pi(\cdot,\cdot,\cdot)$ into this identity proves \eqref{eq:bianchi}.
\end{proof}

We now establish the relationship with Lie algebra homomorphisms. Two constructions will be used. The universal enveloping Lie algebra applies to every Lie triple system, whereas an iterated binary bracket gives a comparison with a prescribed pair of Lie algebras.

For a Lie triple system $L$, let $\mathfrak F(L)$ be the free Lie algebra on its underlying vector space and set
\begin{equation}\label{eq:universal-envelope}
 \mathfrak U(L)=\mathfrak F(L)\big/
 \big\langle [[x,y],z]-[x,y,z]:x,y,z\in L\big\rangle.
\end{equation}
Here $[x,y,z]$ on the right is an element of the generating space $L$. The generators have odd degree in $\mathbb Z_2$. This is an ordinary Lie algebra with a $\mathbb Z_2$-grading. The universal embedding theorem of Jacobson, in the graded formulation of Smirnov, gives the following result \cite{Jacobson,Smirnov}.

\begin{proposition}\label{prop:universal-lie}
The generating map $L\to\mathfrak U(L)$ is injective, and
\begin{equation}\label{eq:universal-grading}
 \mathfrak U(L)_1=L,\qquad
 \mathfrak U(L)_0=[L,L].
\end{equation}
Every morphism $f:\g\to\h$ extends uniquely to a grading-preserving Lie algebra homomorphism
\begin{equation}\label{eq:universal-map}
 \mathfrak U(f):\mathfrak U(\g)\longrightarrow\mathfrak U(\h),
 \qquad \mathfrak U(f)x=fx,\quad
 \mathfrak U(f)[x,y]=[fx,fy].
\end{equation}
Conversely, restriction to the odd component gives a Lie triple system morphism. Thus these two operations give a bijection between the two morphism spaces.
\end{proposition}
\begin{proof}
The injectivity and the decomposition in \eqref{eq:universal-grading} are the universal embedding theorem. For $x,y,z\in\g$, the image of a defining relation under the free Lie homomorphism induced by $f$ is
\[
 [[fx,fy],fz]-f[x,y,z].
\]
In $\mathfrak U(\h)$, this expression is
\[
 [fx,fy,fz]-f[x,y,z],
\]
which vanishes exactly when $f$ preserves the triple product. Thus the free Lie homomorphism descends to the quotient in \eqref{eq:universal-envelope}. It preserves the grading because it takes odd generators to odd generators, and it is unique because those generators generate the source Lie algebra.

Conversely, for a grading-preserving Lie homomorphism from $\mathfrak U(\g)$ to $\mathfrak U(\h)$, preservation of the binary bracket twice gives preservation of $[[x,y],z]$ on odd inputs. Since $[[x,y],z]=[x,y,z]$ in the universal envelope, its odd restriction is a Lie triple morphism. Extension and odd restriction are inverse operations by uniqueness.
\end{proof}

Proposition~\ref{prop:universal-lie} extends every Lie triple morphism through the universal envelopes. For the standard envelope $\InnDer(L)\oplus L$, the assignment $\ad_{x\wedge y}\mapsto\ad_{fx\wedge fy}$ defines a linear map on the even component precisely when it preserves every linear relation among the inner derivations.

Next, let $\mathfrak a$ and $\mathfrak b$ be Lie algebras. Denote their associated Lie triple systems by $\mathfrak a^{\mathrm{tr}}$ and $\mathfrak b^{\mathrm{tr}}$, with products
\[
 \pi(x,y,z)=[[x,y]_{\mathfrak a},z]_{\mathfrak a},\qquad
 \mu(u,v,w)=[[u,v]_{\mathfrak b},w]_{\mathfrak b}.
\]
For a linear map $q:\mathfrak a\to\mathfrak b$, define
\begin{align}
 F(q)(x,y)&=[qx,qy]_{\mathfrak b}-q[x,y]_{\mathfrak a},
       \label{eq:binary-defect}\\
 (\mathcal T_q b)(x,y,z)&=[b(x,y),qz]_{\mathfrak b}
                         +b([x,y]_{\mathfrak a},z),
       \label{eq:binary-triple-map}
\end{align}
where $b:\wedge^2\mathfrak a\to\mathfrak b$. The map $\mathcal T_q$ takes values in $\Hom(\mathfrak a^{\otimes3},\mathfrak b)$.

\begin{proposition}\label{prop:binary-defect}
For every linear map $q:\mathfrak a\to\mathfrak b$, we have
\begin{equation}\label{eq:binary-triple-defect}
 S_{\pi,\mu}(q)=\mathcal T_qF(q).
\end{equation}
Consequently, every Lie algebra homomorphism is a morphism of the associated Lie triple systems. The converse does not hold in general.
\end{proposition}
\begin{proof}
For every $x,y,z\in\mathfrak a$, expanding the binary brackets gives
\begin{equation*}
 [[qx,qy]_{\mathfrak b},qz]_{\mathfrak b} -q[[x,y]_{\mathfrak a},z]_{\mathfrak a} =\bigl[[qx,qy]_{\mathfrak b}-q[x,y]_{\mathfrak a}, qz\bigr]_{\mathfrak b} +[q[x,y]_{\mathfrak a},qz]_{\mathfrak b} -q[[x,y]_{\mathfrak a},z]_{\mathfrak a}.
\end{equation*}
The right-hand side is the full expression obtained from \eqref{eq:binary-triple-map} by substituting
$[q\,\cdot,q\,\cdot]_{\mathfrak b}-q[\cdot,\cdot]_{\mathfrak a}$. This proves \eqref{eq:binary-triple-defect}. In particular, preservation of the binary bracket makes the right-hand side vanish.

On any nonabelian Lie algebra one has
\[
 [[-x,-y],-z]=-[[x,y],z],\qquad
 [-x,-y]=[x,y]\neq-[x,y]
\]
for some $x,y$. Thus $-\id$ preserves the triple product but does not preserve the binary bracket.
\end{proof}

The following condition will give a converse near a fixed homomorphism. Write
\[
 C_{\mathfrak b}(f(\mathfrak a))
 =\{v\in\mathfrak b:[v,fx]_{\mathfrak b}=0
                         \text{ for all }x\in\mathfrak a\}.
\]

\begin{lemma}\label{lem:T-injective}
Let $f:\mathfrak a\to\mathfrak b$ be a Lie algebra homomorphism. If
\begin{equation}\label{eq:Lie-comparison-hyp}
 [\mathfrak a,\mathfrak a]=\mathfrak a,
 \qquad C_{\mathfrak b}(f(\mathfrak a))=0,
\end{equation}
then $\mathcal T_f$ is injective on $\Hom(\wedge^2\mathfrak a,\mathfrak b)$.
\end{lemma}
\begin{proof}
Suppose $b\in\Hom(\wedge^2\mathfrak a,\mathfrak b)$ satisfies
\[
 [b(x,y),fz]_{\mathfrak b}+b([x,y]_{\mathfrak a},z)=0
 \qquad(x,y,z\in\mathfrak a).
\]
If $\sum_i[x_i,y_i]_{\mathfrak a}=0$, summing this equality gives
\[
 \left[\sum_i b(x_i,y_i),fz\right]_{\mathfrak b}=0
 \qquad(z\in\mathfrak a).
\]
The zero-centralizer assumption implies $\sum_i b(x_i,y_i)=0$. Consequently, every linear relation among commutators gives the corresponding relation among the values of $b$.

Applying this observation to the Jacobi identity gives
\begin{equation*}
 b([u,v]_{\mathfrak a},z) =b(u,[v,z]_{\mathfrak a})-b(v,[u,z]_{\mathfrak a}) =-[fu,b(v,z)]_{\mathfrak b}+[fv,b(u,z)]_{\mathfrak b}.
\end{equation*}
Here the second equality follows from the displayed kernel equation and skewness of $b$. Since $\mathfrak a$ is perfect, write any $z\in\mathfrak a$ as $z=\sum_i[x_i,y_i]_{\mathfrak a}$. The same kernel equation gives
\[
 b(u,z)=-\left[fu,\sum_i b(x_i,y_i)\right]_{\mathfrak b}.
\]
Substituting this expression twice in the preceding Jacobi identity yields
\begin{equation*}
 b([u,v]_{\mathfrak a},z) =\left[fu,\left[fv,\sum_i b(x_i,y_i)\right]_{\mathfrak b} \right]_{\mathfrak b} -\left[fv,\left[fu,\sum_i b(x_i,y_i)\right]_{\mathfrak b} \right]_{\mathfrak b} =\left[f[u,v]_{\mathfrak a},\sum_i b(x_i,y_i)\right]_{\mathfrak b}.
\end{equation*}
On the other hand, the formula for $b(u,z)$ with $u$ replaced by $[u,v]_{\mathfrak a}$ gives
\[
 b([u,v]_{\mathfrak a},z)
 =-\left[f[u,v]_{\mathfrak a},\sum_i b(x_i,y_i)\right]_{\mathfrak b}.
\]
Since the characteristic is zero, these two equalities imply
\[
 \left[f[u,v]_{\mathfrak a},\sum_i b(x_i,y_i)\right]_{\mathfrak b}=0.
\]
Perfectness now shows that $\sum_i b(x_i,y_i)$ centralizes $f(\mathfrak a)$, so this sum is zero. The displayed expression for $b(u,z)$ is therefore zero for all $u,z$. Hence $b=0$, proving injectivity.
\end{proof}

\section{Cohomology and formal deformations}

In this section, we define the augmented cohomology of a Lie triple system morphism and study formal deformations. We then compare first cohomology and formal equivalence with those of the universal graded Lie homomorphism. Finally, for morphisms induced by binary Lie brackets, we compare the Chevalley--Eilenberg differentials and obstruction classes.

\begin{definition}[\cite{Yamaguti,XiaShengTang}]
A representation of a Lie triple system $L$ on $V$ is a bilinear map $\theta:L\otimes L\to\End(V)$ such that, with $D(x,y)=\theta(y,x)-\theta(x,y)$,
\begin{align}
 [D(x,y),\theta(u,v)]
   &=\theta([x,y,u],v)+\theta(u,[x,y,v]), \label{eq:rep1}\\
 \theta(x,[u,v,w])
   &=\theta(v,w)\theta(x,u)-\theta(u,w)\theta(x,v)
     +D(u,v)\theta(x,w). \label{eq:rep2}
\end{align}
Equivalently, the following product makes $L\oplus V$ a Lie triple system:
\begin{equation}\label{eq:semidirect}
 [(x,a),(y,b),(z,c)]
 =\bigl([x,y,z],\theta(y,z)a-\theta(x,z)b+D(x,y)c\bigr).
\end{equation}
\end{definition}

\begin{proposition}\label{prop:pullback}
A morphism $f:(\g,\pi)\to(\h,\mu)$ defines a representation on $\h$ by
\begin{equation}\label{eq:pullback}
 \theta_f(x,y)u=\mu(u,fx,fy),\qquad D_f(x,y)u=\mu(fx,fy,u).
\end{equation}
\end{proposition}

\begin{proof}
The adjoint representation of $\h$ is given by $\theta_{\mathrm{ad}}(u,v)w=\mu(w,u,v)$; its representation identities follow from the fundamental and cyclic identities \cite{Yamaguti,XiaShengTang}. For $x,y,u,v\in\g$ and $w\in\h$, substitution in the first adjoint representation identity gives
\begin{align*}
 &\mu(fx,fy,\mu(w,fu,fv))
       -\mu(\mu(fx,fy,w),fu,fv)\\
 &\qquad=\mu(w,\mu(fx,fy,fu),fv)
             +\mu(w,fu,\mu(fx,fy,fv))\\
 &\qquad=\mu(w,f\pi(x,y,u),fv)
             +\mu(w,fu,f\pi(x,y,v)).
\end{align*}
This is \eqref{eq:rep1} for $\theta_f$. Likewise, for $x,u,v,z\in\g$ and $w\in\h$, the second adjoint representation identity gives
\begin{align*}
 \mu(w,fx,f\pi(u,v,z)) &=\mu(w,fx,\mu(fu,fv,fz)) =\mu(\mu(w,fx,fu),fv,fz) -\mu(\mu(w,fx,fv),fu,fz) \\
 &\quad+\mu(fu,fv,\mu(w,fx,fz)).
\end{align*}
This is \eqref{eq:rep2}. Finally, skewness in the first pair and the cyclic identity imply
\[
 \mu(u,fy,fx)-\mu(u,fx,fy)=\mu(fx,fy,u),
\]
which gives the stated formula for $D_f(x,y)u$.
\end{proof}

Put
\begin{equation}\label{eq:aug-cochains}
 C_f^0=\wedge^2\g\oplus\wedge^2\h,
 \qquad C_f^n=\C^{n-1}(\g,\h)\quad(n\geq1).
\end{equation}
Thus $C_f^1=\Hom(\g,\h)$ and $C_f^2$ consists of trilinear maps with the Lie triple symmetries. The numbering in \eqref{eq:aug-cochains} is by cohomological degree. In the traditional odd-arity notation, $C_f^n$ is the Yamaguti cochain space $C_Y^{2n-1}(\g,\h_f)$, where $\h_f$ denotes \eqref{eq:pullback}.

Define
\begin{align}
 (\delta_f^0(X,U))(z)&=\ad_U fz-f\ad_Xz, \label{eq:delta0}\\
 \delta_f^n a&=(-1)^{n-1}\ell_1^f(a)\quad(n\geq1). \label{eq:delta-derived}
\end{align}
To state the full formula, for $X_i=x_i\otimes y_i$ write
\[
 [X_i,X_j]_\pi
   =\pi(x_i,y_i,x_j)\otimes y_j+x_j\otimes\pi(x_i,y_i,y_j).
\]
For $a\in C_f^n$, the differential is
\begin{align}
 &(\delta_f^n a)(X_1,\ldots,X_n,z) =\sum_{1\leq i<j\leq n}(-1)^i a(X_1,\ldots,\widehat X_i,\ldots,X_{j-1},[X_i,X_j]_\pi, X_{j+1},\ldots,X_n,z) \notag \\
 &\quad+\sum_{i=1}^n(-1)^i a(X_1,\ldots,\widehat X_i,\ldots,X_n,\pi(x_i,y_i,z)) \notag \\
 &\quad+\sum_{i=1}^n(-1)^{i+1} \mu\bigl(fx_i,fy_i,a(X_1,\ldots,\widehat X_i,\ldots,X_n,z)\bigr) \notag \\
 &\quad+(-1)^{n+1}\bigl\{ \mu\bigl(a(X_1,\ldots,X_{n-1},x_n),fy_n,fz\bigr) +\mu\bigl(fx_n,a(X_1,\ldots,X_{n-1},y_n),fz\bigr)\bigr\}.\label{eq:full-differential}
\end{align}
The hat indicates omission. Formula \eqref{eq:full-differential} is the Yamaguti differential for \eqref{eq:pullback}. The second endpoint term is obtained from $-\theta_f(x_n,z)a(\ldots,y_n)$ using skewness in the first pair.

In the first two degrees, \eqref{eq:full-differential} reads
\begin{align}
 (\delta_f^1 a)(x,y,z)
 &=\mu(ax,fy,fz)+\mu(fx,ay,fz)+\mu(fx,fy,az)-a\pi(x,y,z),
 \label{eq:d1}\\
 (\delta_f^2 b)(x,y,u,v,w)
 &=\mu(fx,fy,b(u,v,w))-\mu(fu,fv,b(x,y,w)) \notag\\
 &\quad-\mu(b(x,y,u),fv,fw)-\mu(fu,b(x,y,v),fw) \notag\\
 &\quad+b(x,y,\pi(u,v,w))-b(u,v,\pi(x,y,w)) \notag\\
 &\quad-b(\pi(x,y,u),v,w)-b(u,\pi(x,y,v),w).
 \label{eq:d2}
\end{align}

\begin{theorem}\label{thm:complex}
The maps \eqref{eq:delta0}--\eqref{eq:delta-derived} define a cochain complex. In positive degrees this is the Yamaguti complex for $\h_f$, with degrees reindexed as in \eqref{eq:aug-cochains}.
\end{theorem}

\begin{proof}
Insertion in \eqref{eq:delta-derived} gives \eqref{eq:full-differential}, which preserves the cochain symmetries by Proposition~\ref{prop:control}. The unary identity of the twisted $L_\infty[1]$-algebra yields, for $n\geq1$ and $a\in C_f^n$,
\[
 \delta_f^{n+1}(\delta_f^na)
 =(-1)^n(-1)^{n-1}\ell_1^f(\ell_1^f(a))=0.
\]
For degree zero, take $X\in\wedge^2\g$ and $U\in\wedge^2\h$. Since $\ad_U$ is a derivation, the morphism identity gives
\begin{align*}
 (\delta_f^1(\ad_U\circ f))(x,y,z)
 &=\mu(\ad_Ufx,fy,fz)+\mu(fx,\ad_Ufy,fz)\\
 &\quad+\mu(fx,fy,\ad_Ufz)-\ad_Uf\pi(x,y,z)\\
 &=\ad_U\mu(fx,fy,fz)-\ad_Uf\pi(x,y,z)=0.
\end{align*}
Similarly, since $\ad_X$ is a derivation of the source product,
\begin{equation*}
 (\delta_f^1(f\circ\ad_X))(x,y,z) =f\bigl(\pi(\ad_Xx,y,z)+\pi(x,\ad_Xy,z) +\pi(x,y,\ad_Xz)-\ad_X\pi(x,y,z)\bigr)=0.
\end{equation*}
Subtracting these two identities proves
$\delta_f^1(\delta_f^0(X,U))=0$.
\end{proof}

Define
\[
 Z_f^n=\ker\delta_f^n,\qquad
 B_f^n=\im\delta_f^{n-1},\qquad H_f^n=Z_f^n/B_f^n\quad(n\geq1).
\]
We augment the positive-degree Yamaguti complex, which starts with $\Hom(\g,\h)$, by the map $\delta_f^0$. Accordingly,
\begin{equation}\label{eq:degree-comparison}
 H_f^1=Z_Y^1(\g,\h_f)/\im\delta_f^0,
 \qquad H_f^n\cong H_Y^{2n-1}(\g,\h_f)\quad(n\geq2).
\end{equation}
In particular, the second group used below is the odd-arity third Yamaguti cohomology group, whereas the first group has an additional quotient.

\begin{lemma}\label{lem:boundaries}
The first coboundaries can be obtained using target inner derivations alone:
\begin{equation}\label{eq:target-boundaries}
 B_f^1=\{\ad_U\circ f:U\in\wedge^2\h\}.
\end{equation}
\end{lemma}
\begin{proof}
For every $X\in\wedge^2\g$ and $x\in\g$, the morphism identity gives
\[
 f\ad_Xx=\ad_{(\wedge^2f)X}fx.
\]
Consequently, for $U\in\wedge^2\h$,
\[
 (\delta_f^0(X,U))(x)
 =\ad_Ufx-\ad_{(\wedge^2f)X}fx
 =\ad_{U-(\wedge^2f)X}fx.
\]
Thus every first coboundary is obtained by postcomposition with a target inner derivation. Conversely, $\delta_f^0(0,U)=\ad_U\circ f$ for every $U\in\wedge^2\h$, which proves the reverse inclusion.
\end{proof}

There is also a natural pullback of positive cochains from the adjoint complex of $\h$. For $b\in\C^{n-1}(\h,\h)$, set
\[
 (f^*b)(x_1,\ldots,x_{2n-1})=b(fx_1,\ldots,fx_{2n-1}).
\]
Substitution in \eqref{eq:full-differential} gives $\delta_f(f^*b)=f^*(\delta_{\mathrm{ad}}b)$: each source product is taken by $f$ to the corresponding target product, and each action term is \eqref{eq:pullback}. Thus $f^*$ induces maps from adjoint Yamaguti cohomology to $H_f^n$ for $n\geq2$, and from adjoint degree-one cocycles to $H_f^1$ after taking the quotient in \eqref{eq:degree-comparison}.

\begin{definition}
A formal deformation of $f$ with fixed source and target products is a series
\[
 f_t=f+\sum_{r\geq1}t^r f_r\in\Hom(\g,\h)[[t]]
\]
that satisfies the morphism identity over $\K[[t]]$. A deformation of order $N$ satisfies this identity modulo $t^{N+1}$. Two formal deformations are inner equivalent if
\[
 f'_t=B_t f_t A_t^{-1},\qquad A_0=\id_\g,\quad B_0=\id_\h,
\]
where $A_t,B_t$ belong to the $t$-adically completed groups generated by exponentials of positive-order series of inner derivations.
\end{definition}

\begin{proposition}\label{prop:infinitesimal}
The coefficient $f_1$ of a formal deformation lies in $Z_f^1$. Inner equivalent formal deformations have cohomologous first coefficients. Conversely, every element of $Z_f^1$ gives a deformation over $\K[t]/(t^2)$, and its inner equivalence class is determined by its class in $H_f^1$.
\end{proposition}
\begin{proof}
The coefficient of $t$ in the morphism identity for $f_t$ is
\[
 \mu(f_1x,fy,fz)+\mu(fx,f_1y,fz)
       +\mu(fx,fy,f_1z)-f_1\pi(x,y,z)=0.
\]
Hence $f_1\in Z_f^1$. If $f'_t=B_tf_tA_t^{-1}$, taking the coefficient of $t$ gives
\[
 f'_1=f_1+
 \left.\frac{dB_t}{dt}\right|_{t=0}\circ f
 -f\circ\left.\frac{dA_t}{dt}\right|_{t=0}.
\]
The two derivatives are inner derivations by the definition of formal inner equivalence. Formula \eqref{eq:delta0} therefore shows that $f'_1-f_1\in B_f^1$. Conversely, the displayed cocycle equation is precisely the condition that $f+tf_1$ preserve the triple product modulo $t^2$. By Lemma~\ref{lem:boundaries}, every first coboundary has the form $\ad_U\circ f$ for some $U\in\wedge^2\h$. Postcomposition with $\exp(t\ad_U)$ realizes this difference modulo $t^2$, proving the classification.
\end{proof}

\begin{theorem}\label{thm:obstruction}
Let $f_t^{(N)}=\sum_{i=0}^N t^if_i$, with $f_0=f$, be a deformation of order $N\geq1$. Define
\begin{equation}\label{eq:obstruction}
 \mathcal O_{N+1}(x,y,z)
   =\sum_{\substack{i+j+k=N+1\\0\leq i,j,k\leq N}}
       \mu(f_ix,f_jy,f_kz).
\end{equation}
Then $\mathcal O_{N+1}\in Z_f^2$. The deformation extends to order $N+1$ if and only if $[\mathcal O_{N+1}]=0$ in $H_f^2$.
\end{theorem}
\begin{proof}
Since the morphism identity holds through order $N$, one has
\begin{equation*}
 \mu\left(\sum_{i=0}^Nt^if_ix, \sum_{j=0}^Nt^jf_jy, \sum_{k=0}^Nt^kf_kz\right) -\sum_{i=0}^Nt^if_i\pi(x,y,z) =t^{N+1} \sum_{\substack{i+j+k=N+1\\0\leq i,j,k\leq N}} \mu(f_ix,f_jy,f_kz)+O(t^{N+2}).
\end{equation*}
The left-hand side is skew in $x,y$ and has zero cyclic sum in $x,y,z$. Its coefficient of $t^{N+1}$ therefore has the same symmetries and belongs to $C_f^2$.

Apply Lemma~\ref{lem:bianchi} over $\K[[t]]$ to the linear map $\sum_{i=0}^Nt^if_i$. Taking the coefficient of $t^{N+1}$ and expanding \eqref{eq:d2} gives, for $x,y,u,v,w\in\g$,
\begin{align*}
0={}&\sum_{\substack{i+j+k=N+1\\0\leq i,j,k\leq N}}
 \Bigl\{
 \mu(fx,fy,\mu(f_iu,f_jv,f_kw))\\
 &\quad-\mu(fu,fv,\mu(f_ix,f_jy,f_kw))
       -\mu(\mu(f_ix,f_jy,f_ku),fv,fw)\\
 &\quad-\mu(fu,\mu(f_ix,f_jy,f_kv),fw)
       +\mu(f_ix,f_jy,f_k\pi(u,v,w))\\
 &\quad-\mu(f_iu,f_jv,f_k\pi(x,y,w))
       -\mu(f_i\pi(x,y,u),f_jv,f_kw)\\
 &\quad-\mu(f_iu,f_j\pi(x,y,v),f_kw)
 \Bigr\}.
\end{align*}
Thus the trilinear sum in \eqref{eq:obstruction} is a second cocycle. Adding $t^{N+1}f_{N+1}$ to the truncated deformation changes its coefficient of order $N+1$ to
\begin{align*}
 &\sum_{\substack{i+j+k=N+1\\0\leq i,j,k\leq N}}
       \mu(f_ix,f_jy,f_kz)\\
 &\qquad+\mu(f_{N+1}x,fy,fz)+\mu(fx,f_{N+1}y,fz)
       +\mu(fx,fy,f_{N+1}z)-f_{N+1}\pi(x,y,z).
\end{align*}
This expression vanishes for some $f_{N+1}\in\Hom(\g,\h)$ if and only if the trilinear sum belongs to $\im\delta_f^1$, or equivalently its class in $H_f^2$ is zero.
\end{proof}

\begin{corollary}\label{cor:formal-unobstructed}
If $H_f^2=0$, every first cocycle extends to a formal deformation, and every finite-order deformation extends to all orders.
\end{corollary}
\begin{proof}
Every first cocycle gives a deformation modulo $t^2$ by Proposition~\ref{prop:infinitesimal}. At each subsequent order, Theorem~\ref{thm:obstruction} places
\[
 (x,y,z)\longmapsto
 \sum_{\substack{i+j+k=N+1\\0\leq i,j,k\leq N}}
       \mu(f_ix,f_jy,f_kz)
\]
in $\ker\delta_f^2$. Since $H_f^2=0$, this cochain belongs to $\im\delta_f^1$, so a coefficient $f_{N+1}$ can be chosen to cancel it. Induction constructs all coefficients of a formal deformation. The same argument starts from any given finite-order deformation. The resulting series is defined in $\K[[t]]$.
\end{proof}

\begin{theorem}\label{thm:formal-rigid}
If $H_f^1=0$, every formal deformation of $f$ is inner equivalent to the constant deformation.
\end{theorem}
\begin{proof}
Suppose the coefficients of orders $1,\ldots,r-1$ have already been removed. The coefficient of $t^r$ in the morphism identity is then
\[
 \mu(f_rx,fy,fz)+\mu(fx,f_ry,fz)
       +\mu(fx,fy,f_rz)-f_r\pi(x,y,z)=0.
\]
Since $H_f^1=0$, Lemma~\ref{lem:boundaries} gives an element $U_r\in\wedge^2\h$ such that
\[
 f_rx=\ad_{U_r}fx\qquad(x\in\g).
\]
Postcomposition with $\exp(-t^r\ad_{U_r})$ changes the coefficient of order $r$ to $f_r-\ad_{U_r}\circ f=0$ and leaves all lower coefficients unchanged. Repeating this procedure produces a product of formal inner automorphisms. It is well defined in the $t$-adic topology, since only finitely many factors affect a prescribed coefficient. After all orders have been removed, the transformed formal deformation is the constant map $f$.
\end{proof}

\begin{example}\label{ex:cubic}
Let $\g=\R^2$ have zero triple product. Give $\h=\R^2$ its Euclidean product
\begin{equation}\label{eq:euclidean-lts}
 \mu(u,v,w)=\langle v,w\rangle u-\langle u,w\rangle v.
\end{equation}
This is a Lie triple system: it has the required skew and cyclic symmetries, and each $\ad_{u\wedge v}$ is skew-adjoint, hence preserves \eqref{eq:euclidean-lts}. Take $f=0$ and identify the two underlying vector spaces. Then $\delta_f^1=\delta_f^2=0$, because both the source product and the coefficient action vanish. Also $\delta_f^0=0$.

The identity map $a:\g\to\h$ is therefore a first cocycle, and $ta$ solves the morphism equation modulo $t^3$. Its defect is
\[
 S_{0,\mu}(ta)(x,y,z)=t^3\mu(ax,ay,az).
\]
For an orthonormal basis $e_1,e_2$, its third-order coefficient at $(e_1,e_2,e_1)$ is $-e_2$. Since $\delta_f^1=0$, no coefficient of order three can remove it. Even adding an arbitrary term $t^2b$ leaves this third-order coefficient unchanged. Thus $a$ is not the infinitesimal of a formal deformation.

Here $C_f^2$ has dimension four: a cochain is determined by its two values on $(e_1,e_2,e_1)$ and $(e_1,e_2,e_2)$, and the cyclic condition is automatic in dimension two. Hence $H_f^2=C_f^2\neq0$. The cubic term in \eqref{eq:twisted-MC} produces this nonzero third-order obstruction.
\end{example}

We now compare this cohomology and the formal deformation problem with those for Lie algebra homomorphisms. First consider the universal envelopes. Put $\Phi=\mathfrak U(f)$ and let
\[
 C_{\mathrm{CE},\mathrm{gr}}^n(\Phi)
 =\Hom_{\mathrm{gr}}\bigl(\wedge^n\mathfrak U(\g),\mathfrak U(\h)\bigr),
 \qquad C_{\mathrm{CE},\mathrm{gr}}^0(\Phi)=\mathfrak U(\h)_0.
\]
The subscript $\mathrm{gr}$ means that a cochain preserves total parity. The coefficient action is $v\cdot w=[\Phi v,w]$. The Chevalley--Eilenberg differential preserves these cochains. In degrees zero and one it is
\begin{equation}\label{eq:graded-CE-low}
 (d_\Phi v)(X)=[\Phi X,v],\qquad (d_\Phi A)(X,Y)=[\Phi X,A(Y)]+[A(X),\Phi Y]-A([X,Y]).
\end{equation}
Denote its cohomology by $H_{\mathrm{CE},\mathrm{gr}}^n(\Phi)$.

\begin{theorem}\label{thm:graded-deformation-comparison}
Restriction to the odd component induces an isomorphism
\begin{equation}\label{eq:graded-first-iso}
 H_{\mathrm{CE},\mathrm{gr}}^1(\mathfrak U(f))\cong H_f^1.
\end{equation}
Moreover, $f_t\mapsto\mathfrak U(f_t)$ gives a bijection between formal deformations of $f$ and grading-preserving formal deformations of $\mathfrak U(f)$. This bijection holds at every finite order and preserves inner equivalence, where the Lie algebra equivalence is generated by even inner derivations. In particular, formal rigidity of $f$ is equivalent to formal rigidity of $\mathfrak U(f)$ in this graded sense.
\end{theorem}
\begin{proof}
Let $A$ be a grading-preserving first Chevalley--Eilenberg cocycle for $\mathfrak U(f)$. The cocycle identity means that $\mathfrak U(f)+\varepsilon A$ preserves the binary bracket over $\K[\varepsilon]/(\varepsilon^2)$. Its odd restriction therefore preserves the triple product. The coefficient of $\varepsilon$ on odd inputs is
\[
 \mu(Ax,fy,fz)+\mu(fx,Ay,fz)
       +\mu(fx,fy,Az)-A\pi(x,y,z)=0.
\]
Thus the odd restriction of $A$ belongs to $\ker\delta_f^1$.

Conversely, for $a\in\ker\delta_f^1$, the map $f+\varepsilon a$ preserves the defining triple relations modulo $\varepsilon^2$. By the universal presentation, it has a unique graded Lie extension, whose values on odd generators and their brackets are
\begin{equation}\label{eq:infinitesimal-envelope}
 \mathfrak U(f+\varepsilon a)(x) =fx+\varepsilon a(x),\qquad \mathfrak U(f+\varepsilon a)([x,y]) =[fx,fy]+\varepsilon\bigl([a(x),fy]+[fx,a(y)]\bigr).
\end{equation}
Its coefficient of $\varepsilon$ is the unique graded first cocycle with odd restriction $a$. Hence odd restriction identifies the first cocycle spaces.

Every even element of $\mathfrak U(\h)$ is a finite sum of brackets of odd elements. For such a sum, its inner action on the odd component is
\[
 \left[\sum_i[u_i,v_i],fx\right]
 =\sum_i[[u_i,v_i],fx]
 =\sum_i[u_i,v_i,fx].
\]
Lemma~\ref{lem:boundaries} identifies these maps with $B_f^1$. The differential in \eqref{eq:graded-CE-low} gives the negatives of these inner maps, so its restricted image is also $B_f^1$. Thus restriction also identifies the first boundary spaces and proves \eqref{eq:graded-first-iso}.

The same presentation argument applies over $\K[t]/(t^{N+1})$ at every finite order. Passing coefficientwise to the inverse limit gives the formal correspondence, with $t$-adic completions in infinite dimension. On the odd component, the even Lie inner derivation $\ad_{[u,v]}$ acts by $x\mapsto[u,v,x]$, which is exactly $\ad_{u\wedge v}$. Conversely, every Lie triple inner derivation is a finite sum of such actions. Exponentiating positive-order series therefore gives the same formal inner equivalences on the source and target odd components. This proves all assertions.
\end{proof}

Consequently, $H_f^1=0$ implies graded formal rigidity of $\mathfrak U(f)$. If $H_f^2=0$, every graded Lie infinitesimal and every finite-order graded deformation of $\mathfrak U(f)$ extends to a formal deformation, by Corollary~\ref{cor:formal-unobstructed} and Theorem~\ref{thm:graded-deformation-comparison}. These graded extension results follow from the Lie triple obstruction criterion in $H_f^2$.

For the second comparison, let $f:\mathfrak a\to\mathfrak b$ be a Lie algebra homomorphism, and let $d_f$ be the Chevalley--Eilenberg differential for the action $x\cdot u=[fx,u]_{\mathfrak b}$. For $a:\mathfrak a\to\mathfrak b$, we have
\begin{equation}\label{eq:ordinary-CE-first}
 (d_fa)(x,y)=[fx,ay]_{\mathfrak b}+[ax,fy]_{\mathfrak b}
                    -a[x,y]_{\mathfrak a}.
\end{equation}
In the following statements, $H_f^n$ denotes the cohomology of the associated Lie triple system morphism.

\begin{proposition}\label{prop:ordinary-cohomology}
For every $a:\mathfrak a\to\mathfrak b$, we have
\begin{equation}\label{eq:differential-comparison}
 \delta_f^1a=\mathcal T_f(d_fa).
\end{equation}
Hence $Z_{\mathrm{CE}}^1(f)\subseteq Z_f^1$. The map $b\mapsto\mathcal T_fb$ induces a linear map
\begin{equation}\label{eq:second-comparison}
 H_{\mathrm{CE}}^2(f)\longrightarrow H_f^2.
\end{equation}
If Eq.~\eqref{eq:Lie-comparison-hyp} holds, then the first cocycle spaces are equal and the map in Eq.~\eqref{eq:second-comparison} is injective. If, in addition, $\mathfrak b$ is perfect, then
\begin{equation}\label{eq:ordinary-first-iso}
 H_{\mathrm{CE}}^1(f)\cong H_f^1.
\end{equation}
For every Lie algebra homomorphism $f$, the first boundary spaces are
\begin{equation}\label{eq:ordinary-boundaries}
 B_{\mathrm{CE}}^1(f)=\{[v,f(\cdot)]:v\in\mathfrak b\},\qquad
 B_f^1=\{[v,f(\cdot)]:v\in[\mathfrak b,\mathfrak b]\}.
\end{equation}
\end{proposition}
\begin{proof}
For $a:\mathfrak a\to\mathfrak b$, direct expansion of the first Lie triple differential gives
\begin{align*}
 (\delta_f^1a)(x,y,z) &=[[ax,fy]_{\mathfrak b},fz]_{\mathfrak b} +[[fx,ay]_{\mathfrak b},fz]_{\mathfrak b} +[[fx,fy]_{\mathfrak b},az]_{\mathfrak b} -a[[x,y]_{\mathfrak a},z]_{\mathfrak a} \\
 &=\bigl[[fx,ay]_{\mathfrak b}+[ax,fy]_{\mathfrak b} -a[x,y]_{\mathfrak a},fz\bigr]_{\mathfrak b} +[f[x,y]_{\mathfrak a},az]_{\mathfrak b} +[a[x,y]_{\mathfrak a},fz]_{\mathfrak b} -a[[x,y]_{\mathfrak a},z]_{\mathfrak a}.
\end{align*}
This is the full expression in \eqref{eq:binary-triple-map} applied to the binary differential in \eqref{eq:ordinary-CE-first}, proving \eqref{eq:differential-comparison}.

Let $b$ be a Chevalley--Eilenberg second cocycle. On $\mathfrak a\oplus\mathfrak b$, regard the second summand as an abelian vector space and use the bracket
\[
 [(x,u),(y,v)]
 =\bigl([x,y]_{\mathfrak a},
       [fx,v]_{\mathfrak b}-[fy,u]_{\mathfrak b}+b(x,y)\bigr).
\]
The representation identity and the second cocycle equation give the Jacobi identity for this bracket. Its iterated bracket is
\begin{align*}
 &[(x,u),(y,v),(z,w)] =\Bigl([[x,y]_{\mathfrak a},z]_{\mathfrak a}, [[u,fy]_{\mathfrak b},fz]_{\mathfrak b} +[[fx,v]_{\mathfrak b},fz]_{\mathfrak b} \\
 &\hspace{5em}+[[fx,fy]_{\mathfrak b},w]_{\mathfrak b} +[b(x,y),fz]_{\mathfrak b}+b([x,y]_{\mathfrak a},z)\Bigr).
\end{align*}
The skew and cyclic identities show that the trilinear map
\[
 (x,y,z)\longmapsto[b(x,y),fz]_{\mathfrak b}
                         +b([x,y]_{\mathfrak a},z)
\]
belongs to $C_f^2$. Applying the fundamental identity to inputs in the first summand, and expanding the second component, gives
\[
 \delta_f^2\bigl([b(\cdot,\cdot),f\,\cdot]_{\mathfrak b}
             +b([\cdot,\cdot]_{\mathfrak a},\cdot)\bigr)=0.
\]
The terms are those of \eqref{eq:d2} with the displayed trilinear expression substituted throughout. Therefore second Lie cocycles give second Lie triple cocycles. The first differential identity already proved shows that second coboundaries also map to second coboundaries. This proves \eqref{eq:second-comparison}.

Under \eqref{eq:Lie-comparison-hyp}, Lemma~\ref{lem:T-injective} and \eqref{eq:differential-comparison} give equality of the first cocycle spaces. If the class of a second Lie cocycle maps to zero, there is $a:\mathfrak a\to\mathfrak b$ such that
\[
 \mathcal T_f\bigl(b(\cdot,\cdot)
       -[f\,\cdot,a\,\cdot]_{\mathfrak b}
       -[a\,\cdot,f\,\cdot]_{\mathfrak b}
       +a[\cdot,\cdot]_{\mathfrak a}\bigr)=0.
\]
Injectivity implies
\[
 b(x,y)=[fx,ay]_{\mathfrak b}+[ax,fy]_{\mathfrak b}
                                -a[x,y]_{\mathfrak a},
\]
so the original class is zero. This proves injectivity in degree two. Finally, a Lie triple inner action on the target is $x\mapsto[[u,v]_{\mathfrak b},fx]_{\mathfrak b}$, whereas a Lie inner action is $x\mapsto[v,fx]_{\mathfrak b}$. Their spans are exactly the two spaces in \eqref{eq:ordinary-boundaries}. If $\mathfrak b$ is perfect, these boundary spaces agree. Together with equality of the first cocycle spaces, this proves \eqref{eq:ordinary-first-iso}.
\end{proof}

For a finite-order Lie homomorphism deformation $f_t^{(N)}=\sum_{i=0}^N t^if_i$, its binary obstruction is
\[
 \Omega_{N+1}(x,y)
 =\sum_{\substack{i+j=N+1\\0\leq i,j\leq N}}
                         [f_ix,f_jy]_{\mathfrak b}.
\]
For an arbitrary linear map $q$, let $d_q$ be the linear operator given by the usual degree-two Chevalley--Eilenberg formula with $x\cdot u=[qx,u]_{\mathfrak b}$. The Jacobi identities of the two binary brackets give $d_qF(q)=0$. Comparing the coefficient of $t^{N+1}$ yields $d_f\Omega_{N+1}=0$. Adding a term $t^{N+1}a$ changes the binary obstruction to $\Omega_{N+1}+d_fa$. On the other hand, Eq.~\eqref{eq:binary-triple-defect} gives
\begin{equation}\label{eq:obstruction-comparison}
 \mathcal O_{N+1}=\mathcal T_f\Omega_{N+1}.
\end{equation}
Indeed, $F(f_t^{(N)})$ has no terms of orders at most $N$, so only the constant term $\mathcal T_f$ contributes at order $N+1$. Thus Eq.~\eqref{eq:second-comparison} takes the Lie obstruction class to the Lie triple obstruction class. Under Eq.~\eqref{eq:Lie-comparison-hyp}, a Lie finite-order deformation extends by one order in the Lie category if and only if it extends in the Lie triple category. Applying the same argument to the first nonzero coefficient of a binary defect shows that every formal Lie triple deformation based at $f$ is then a Lie homomorphism deformation. If both Lie algebras are perfect, their Lie and Lie triple inner equivalences also agree, by Eq.~\eqref{eq:ordinary-boundaries} applied to source and target. Formal inner rigidity is therefore equivalent in the two categories under these hypotheses.

\section{Local rigidity and stability of morphisms}

In this section, $\g$ and $\h$ are finite-dimensional real vector spaces. We study smooth deformations, local rigidity and stability by means of the augmented cohomology. Subsequently, we obtain the corresponding results for universal graded Lie homomorphisms and establish the local comparison for triple products induced by binary Lie brackets.

\begin{definition}
A smooth deformation of $f$ is a smooth curve $f_t$ of morphisms with $f_0=f$. Two such curves are inner equivalent if $f'_t=B_tf_tA_t^{-1}$ for smooth families of inner automorphisms with $A_0=B_0=\id$. A morphism is locally inner rigid if every sufficiently close morphism belongs to its inner orbit.
\end{definition}

Here $\Inn(L)$ means the group generated by $\exp(D)$, $D\in\InnDer(L)$. These are automorphisms because $D$ is a derivation. One may equivalently use the action of the simply connected Lie group integrating $\InnDer(L)$. The argument below uses its local exponential parametrization.

\begin{proposition}\label{prop:smooth-inf}
For a smooth deformation $f_t$, the derivative $\dot f_0$ lies in $Z_f^1$. Inner equivalent deformations have the same class in $H_f^1$. The tangent space to the inner orbit through $f$ is $B_f^1$.
\end{proposition}
\begin{proof}
Differentiating the morphism identity at $t=0$ gives
\[
 \mu(\dot f_0x,fy,fz)+\mu(fx,\dot f_0y,fz)
       +\mu(fx,fy,\dot f_0z)-\dot f_0\pi(x,y,z)=0.
\]
Thus $\dot f_0\in\ker\delta_f^1$. Differentiating $f'_t=B_tf_tA_t^{-1}$ at zero gives
\[
 \dot f'_0-\dot f_0
 =\left.\frac{dB_t}{dt}\right|_{t=0}\circ f
       -f\circ\left.\frac{dA_t}{dt}\right|_{t=0}
 \in\im\delta_f^0.
\]
Conversely, for $X\in\wedge^2\g$ and $U\in\wedge^2\h$, the curve
\[
 t\longmapsto\exp(t\ad_U)\circ f\circ\exp(-t\ad_X)
\]
has derivative $\ad_U\circ f-f\circ\ad_X$. Its derivatives exhaust the first boundary space by \eqref{eq:delta0}, proving the tangent-space assertion. Lemma~\ref{lem:boundaries} shows that taking $X=0$ already gives all these tangent vectors.
\end{proof}

\begin{theorem}\label{thm:local-rigid}
If $H_f^1=0$, there is an open neighborhood $U$ of $f$ in $\Hom(\g,\h)$ and a smooth map
\[
 D:U\longrightarrow\InnDer(\h),\qquad D(f)=0,
\]
such that every morphism $q\in U$ satisfies
\begin{equation}\label{eq:rigid-normalization}
 q=\exp(D(q))f.
\end{equation}
In particular, $f$ is locally inner rigid, and every sufficiently small smooth deformation is smoothly inner equivalent to the constant one.
\end{theorem}
\begin{proof}
Choose a vector-space complement $N$ to $\im\delta_f^0$ in $\Hom(\g,\h)$, and choose a complement $\mathfrak a$ to
\[
 \{D\in\InnDer(\h):D\circ f=0\}
\]
in $\InnDer(\h)$. Lemma~\ref{lem:boundaries} implies that $D\mapsto D\circ f$ restricts to an isomorphism from $\mathfrak a$ to $\im\delta_f^0$. The smooth map
\[
 \mathfrak a\times N\longrightarrow\Hom(\g,\h),
 \qquad(D,n)\longmapsto\exp(D)\circ(f+n)
\]
has derivative $(D,n)\mapsto D\circ f+n$ at $(0,0)$, which is an isomorphism. The inverse function theorem therefore gives smooth coordinates $D(q)$ and $n(q)$ satisfying
\[
 q=\exp(D(q))\circ(f+n(q)),\qquad D(f)=0,\quad n(f)=0.
\]
For every target automorphism $B$ and every linear map $q$, one has
\begin{equation*}
 \mu(Bqx,Bqy,Bqz)-Bq\pi(x,y,z) =B\bigl(\mu(qx,qy,qz)-q\pi(x,y,z)\bigr).
\end{equation*}
Hence $f+n(q)$ is a morphism whenever $q$ is a morphism. Differentiating the defining polynomial gives
\begin{align}
 (dS_{\pi,\mu})_f(a)(x,y,z)
 &=\mu(ax,fy,fz)+\mu(fx,ay,fz)\notag\\
 &\quad+\mu(fx,fy,az)-a\pi(x,y,z)
 =(\delta_f^1a)(x,y,z).
 \label{eq:linearization}
\end{align}
The hypothesis $H_f^1=0$ gives $\ker\delta_f^1=\im\delta_f^0$. The displayed linear map is therefore injective on $N$. In finite dimension, there is $c>0$ such that, for $n\in N$,
\begin{equation*}
 \bigl\|\mu(n\,\cdot,f\,\cdot,f\,\cdot) +\mu(f\,\cdot,n\,\cdot,f\,\cdot) +\mu(f\,\cdot,f\,\cdot,n\,\cdot) -n\pi(\cdot,\cdot,\cdot)\bigr\| \geq c\|n\|.
\end{equation*}
Taylor expansion consequently gives, for a constant $C>0$ and sufficiently small $n\in N$,
\begin{equation*}
 \bigl\|\mu((f+n)\,\cdot,(f+n)\,\cdot,(f+n)\,\cdot) -(f+n)\pi(\cdot,\cdot,\cdot)\bigr\| \geq c\|n\|-C\|n\|^2.
\end{equation*}
Thus $f+n$ can be a morphism on this slice only when $n=0$. It follows that $n(q)=0$ for every nearby morphism, so $q=\exp(D(q))\circ f$. Composing $D$ with a sufficiently small smooth deformation and postcomposing by $\exp(-D(f_t))$ gives its smooth inner equivalence with the constant deformation.
\end{proof}

The next elementary lemma is a local form of the exactness argument used in \cite{CrainicSchaetzStruchiner}. We include its proof to make the parameter dependence explicit.

\begin{lemma}\label{lem:exact-zero}
Let $V,E,F$ be finite-dimensional real vector spaces, let $s$ be a smooth map near $x_0\in V$ with values in $E$, and let $B(x):E\to F$ depend smoothly on $x$. Suppose
\[
 s(x_0)=0,\qquad B(x)s(x)=0,
 \qquad \im(ds)_{x_0}=\ker B(x_0).
\]
Then $s^{-1}(0)$ is a smooth submanifold near $x_0$, with tangent space $\ker(ds)_{x_0}$.

Suppose, in addition, that $s_p$ and $B_p$ depend smoothly on a parameter $p$ in an ambient finite-dimensional open set, with $s_{p_0}=s$ and $B_{p_0}=B$. Assume that $B_p(x)s_p(x)=0$ for parameters in a subset $\mathcal P$ containing $p_0$. For every $p\in\mathcal P$ close to $p_0$, there is a zero $x(p)$ close to $x_0$. The choice $x(p)$ is the restriction of a smooth map on the ambient parameter space.
\end{lemma}
\begin{proof}
Choose a vector-space complement $C$ with
\[
 E=\im(ds)_{x_0}\oplus C.
\]
Since $\ker B(x_0)=\im(ds)_{x_0}$, the restriction $B(x_0)|_C$ is injective. This injectivity persists for $B_p(x)|_C$ when $p$ and $x$ are sufficiently close to their base values. Choose a vector-space complement $W$ with
\[
 V=\ker(ds)_{x_0}\oplus W.
\]
Then $(ds)_{x_0}|_W:W\to\im(ds)_{x_0}$ is an isomorphism. Projecting $E$ onto $\im(ds)_{x_0}$ along $C$, the equation
\[
 \pr_{\im(ds)_{x_0}}s_p(x_0+k+w)=0,
 \qquad k\in\ker(ds)_{x_0},\quad w\in W,
\]
has a unique solution $w=w(p,k)$ near zero by the inverse function theorem. This solution is smooth in $(p,k)$ and satisfies $w(p_0,0)=0$. The projected equation implies
\[
 s_p(x_0+k+w(p,k))\in C.
\]
For $p\in\mathcal P$, the compatibility identity also gives
\[
 B_p(x_0+k+w(p,k))s_p(x_0+k+w(p,k))=0.
\]
Injectivity of $B_p(x_0+k+w(p,k))|_C$ therefore implies $s_p(x_0+k+w(p,k))=0$.

When $p=p_0$, every nearby zero of $s$ satisfies the projected equation. Hence the zero set is exactly the graph of $k\mapsto w(p_0,k)$ over $\ker(ds)_{x_0}$. Differentiating the projected equation in a direction $k\in\ker(ds)_{x_0}$ gives
\[
 \pr_{\im(ds)_{x_0}}(ds)_{x_0}
       \bigl(k+(d_kw)_{(p_0,0)}k\bigr)=0.
\]
The derivative vanishes on $k$, and its restriction to $W$ is an isomorphism, so $(d_kw)_{(p_0,0)}k=0$. This proves the tangent-space assertion. Finally, $p\mapsto x_0+w(p,0)$ is smooth on the ambient parameter neighborhood and is a zero of $s_p$ whenever $p\in\mathcal P$.
\end{proof}

Let $\operatorname{Mor}(\g,\h)$ denote the zero set of \eqref{eq:defect} for the fixed products $\pi,\mu$.

\begin{theorem}\label{thm:smoothness}
If $H_f^2=0$, then $\operatorname{Mor}(\g,\h)$ is a smooth submanifold near $f$, with
\[
 T_f\operatorname{Mor}(\g,\h)=Z_f^1,
 \qquad \dim_f\operatorname{Mor}(\g,\h)=\dim Z_f^1.
\]
Every $a\in Z_f^1$ is the derivative of a smooth deformation of $f$.
\end{theorem}
\begin{proof}
Apply Lemma~\ref{lem:exact-zero} to the polynomial map
\begin{equation*}
 \Hom(\g,\h)\longrightarrow\C^1(\g,\h),\qquad q\longmapsto \bigl((x,y,z)\longmapsto\mu(qx,qy,qz)-q\pi(x,y,z)\bigr),
\end{equation*}
and to the family of linear maps from $\C^1(\g,\h)$ to $\C^2(\g,\h)$ given by
\[
 b\longmapsto\mathsf P[\pi,\widehat b]_{\B}
       +\tfrac12\mathsf P[[[\mu,\widehat q]_{\B},\widehat q]_{\B},
                              \widehat b]_{\B}.
\]
The compatibility required by the lemma is exactly Lemma~\ref{lem:bianchi}. At $q=f$, the derivative of the polynomial map sends $a$ to
\[
 (x,y,z)\longmapsto
 \mu(ax,fy,fz)+\mu(fx,ay,fz)
       +\mu(fx,fy,az)-a\pi(x,y,z),
\]
which is $\delta_f^1a$. By \eqref{eq:delta-derived}, the displayed family of linear maps at $q=f$ sends $b$ to $-\delta_f^2b$. Thus the necessary exactness is
\[
 \im\delta_f^1=\ker\delta_f^2,
\]
which follows from $H_f^2=0$. Lemma~\ref{lem:exact-zero} proves local smoothness and identifies the tangent space with $\ker\delta_f^1=Z_f^1$. Every vector in the tangent space of a smooth submanifold is the derivative of a smooth curve through the base point, proving the final assertion.
\end{proof}

\begin{definition}
A morphism $f:(\g,\pi)\to(\h,\mu)$ is stable under simultaneous perturbations if, for every neighborhood $U$ of $f$, there is a neighborhood $W$ of $(\pi,\mu)$ such that every pair of Lie triple system products $(\pi',\mu')\in W$ admits a morphism $f'\in U$ from $(\g,\pi')$ to $(\h,\mu')$.
\end{definition}

\begin{theorem}\label{thm:stability}
If $H_f^2=0$, then $f$ is stable under simultaneous perturbations. More precisely, a nearby morphism can be chosen as the restriction of a smooth function of the two product tensors. Consequently, for any smooth families $\pi_t,\mu_t$ of Lie triple system products with $(\pi_0,\mu_0)=(\pi,\mu)$, there is, after restricting the parameter interval, a smooth family $f_t$ of morphisms with $f_0=f$.
\end{theorem}
\begin{proof}
Use the ambient parameter vector space
\[
 \C^1(\g,\g)\oplus\C^1(\h,\h).
\]
For arbitrary $(\pi',\mu')$ in this space, consider the polynomial map
\[
 q\longmapsto
 \bigl((x,y,z)\longmapsto\mu'(qx,qy,qz)-q\pi'(x,y,z)\bigr)
\]
and the family of linear maps
\[
 b\longmapsto\mathsf P[\pi',\widehat b]_{\B}
       +\tfrac12\mathsf P[[[\mu',\widehat q]_{\B},\widehat q]_{\B},
                              \widehat b]_{\B}.
\]
Both depend polynomially on $(\pi',\mu',q)$. Whenever $\pi'$ and $\mu'$ satisfy the Lie triple identities, Lemma~\ref{lem:bianchi} gives the compatibility identity for these two maps. At $(\pi,\mu,f)$, their linearization and value are $\delta_f^1$ and $-\delta_f^2$, respectively, as shown in Theorem~\ref{thm:smoothness}. The condition $H_f^2=0$ therefore gives the exactness required by Lemma~\ref{lem:exact-zero}.

The parameter part of that lemma gives a smooth choice $q(\pi',\mu')$, with $q(\pi,\mu)=f$, such that
\begin{equation*}
 \mu'\bigl(q(\pi',\mu')x,q(\pi',\mu')y,q(\pi',\mu')z\bigr) =q(\pi',\mu')\pi'(x,y,z)
\end{equation*}
whenever the two parameter tensors are Lie triple products. Continuity of $q$ gives the neighborhood condition in the definition of stability. For smooth families $\pi_t,\mu_t$, the composition $t\mapsto q(\pi_t,\mu_t)$ gives the required smooth family of morphisms through $f$.
\end{proof}

\begin{corollary}
Under the hypothesis of Theorem~\ref{thm:stability}, $f$ is stable when only the target product, or only the source product, is perturbed. If $f$ is injective, surjective or invertible, respectively, the nearby morphisms can be chosen with the same property.
\end{corollary}
\begin{proof}
For perturbations of the target product alone, apply Theorem~\ref{thm:stability} to pairs $(\pi,\mu')$. For perturbations of the source product alone, apply it to pairs $(\pi',\mu)$. In finite dimension, the set of injective linear maps is open, as is the set of surjective linear maps. When the dimensions agree, the set of invertible linear maps is open. The nearby morphisms supplied by the theorem converge to $f$ and therefore retain each of these properties whenever $f$ has it.
\end{proof}

\begin{example}\label{ex:rigid}
Let $T=\R^2$ with \eqref{eq:euclidean-lts}, and take $f=\id_T$. In an orthonormal basis put
\[
 J=\begin{pmatrix}0&1\\-1&0\end{pmatrix},\qquad
 [x,y,z]=\omega(x,y)Jz,
\]
where $\omega(e_1,e_2)=1$. For an endomorphism
$A=\left(\begin{smallmatrix}a&b\\c&d\end{smallmatrix}\right)$,
\begin{equation}\label{eq:example-d1}
 (\delta_f^1A)(e_1,e_2,z)
   =\bigl((\tr A)J+JA-AJ\bigr)z
   =\begin{pmatrix}b+c&2d\\-2a&-b-c\end{pmatrix}z.
\end{equation}
Thus $Z_f^1=\R J$. Since $\InnDer(T)=\R J$, we have $B_f^1=\R J$ and $H_f^1=0$.

Every $b\in C_f^2$ has the form $b(x,y,z)=\omega(x,y)Kz$ for a unique $K\in M_2(\R)$. Formula \eqref{eq:d2} gives
\begin{equation}\label{eq:example-d2}
 (\delta_f^2b)(x,y,u,v,w)
   =-\omega(x,y)\omega(u,v)(\tr K)Jw.
\end{equation}
Consequently, $Z_f^2$ is the three-dimensional space $\tr K=0$. By \eqref{eq:example-d1}, this is exactly $B_f^2$, so $H_f^2=0$. The identity morphism is locally inner rigid and stable, and its nearby morphism space has dimension one. The tangent direction is $J$, corresponding to the rotation family $\exp(tJ)$.
\end{example}

\begin{example}\label{ex:line}
Let $\g=\R e$ have zero product and let $\h$ be any finite-dimensional real Lie triple system. Every linear map $f:\g\to\h$ is a morphism, because $[u,u,u]=0$. Also $C_f^2=0$, since a trilinear cochain is skew in its first two entries. Hence $H_f^2=0$ and every such morphism is stable.

For $\h=T$ from Example~\ref{ex:rigid} and $f(e)=e_1$, one has
\[
 Z_f^1\cong\R^2,\qquad B_f^1=\R e_2,\qquad \dim H_f^1=1.
\]
The family $f_t(e)=(1+t)e_1$ cannot be obtained from $f$ by inner automorphisms: the inner automorphisms of $T$ are rotations and preserve length, whereas the source has no nonzero inner derivations. Thus stability does not imply rigidity.
\end{example}

We next apply these results to the universal enveloping Lie algebra homomorphism. In the following statements, a graded deformation preserves the fixed $\mathbb Z_2$-decomposition. Graded inner equivalence means equivalence generated by exponentials of even Lie inner derivations.

\begin{theorem}\label{thm:graded-local-comparison}
Let $\g$ and $\h$ be finite-dimensional real Lie triple systems, and let $\Phi=\mathfrak U(f)$. Restriction to the odd component identifies the local morphism space at $\Phi$ in the graded Lie category with the local Lie triple morphism space at $f$. Both the map and its inverse extend to polynomial maps on the ambient spaces of linear maps. They preserve smooth deformations and the corresponding inner orbits. Consequently, $f$ is locally inner rigid if and only if $\Phi$ is locally inner rigid in the graded sense.

If $H_f^1=0$, then $\Phi$ is locally graded inner rigid. If $H_f^2=0$, then its graded morphism space is locally smooth, and every graded Lie first cocycle is the derivative of a smooth graded deformation.
\end{theorem}
\begin{proof}
The even component $\mathfrak U(\g)_0$ is a quotient of $\wedge^2\g$, and the odd component is $\g$. Thus the two universal envelopes are finite-dimensional. Choose pairs $x_i,y_i\in\g$ such that the brackets $[x_i,y_i]$ form a basis of $\mathfrak U(\g)_0$. For an arbitrary linear map $q:\g\to\h$, the assignments
\[
 x\longmapsto qx\quad(x\in\g),\qquad
 [x_i,y_i]\longmapsto[qx_i,qy_i]
\]
on the odd space and the chosen even basis define a grading-preserving linear map on the source envelope. Its coefficients depend polynomially on the coefficients of $q$. On Lie triple morphisms it is precisely the unique universal extension from Proposition~\ref{prop:universal-lie}. Conversely, odd restriction is a linear map on the ambient space of graded linear maps. The two operations are inverse on the corresponding morphism loci.

These polynomial maps identify smooth curves on the two loci. Whenever one locus is a smooth submanifold, their compositions restrict to the identity and their differentials are inverse on the corresponding tangent spaces; the same local parametrizations consequently give the submanifold structure on the other locus.

For even inner derivations, the restriction formula is
\[
 \left.\ad_{\sum_i[u_i,v_i]}\right|_{\h}(x)
 =\sum_i[u_i,v_i,x]
 =\left(\sum_i\ad_{u_i\wedge v_i}\right)(x).
\]
Every Lie triple inner derivation arises in this way. Exponentiation therefore identifies the inner orbits on the two morphism spaces, proving the rigidity equivalence. The remaining assertions follow from Theorems~\ref{thm:local-rigid}, \ref{thm:smoothness} and~\ref{thm:graded-deformation-comparison}.
\end{proof}

For stability, the products on the envelopes are also allowed to vary. A graded Lie homomorphism is called graded stable if it persists under every sufficiently small pair of Lie bracket perturbations preserving the fixed grading on its source and target.

\begin{theorem}\label{thm:graded-stability}
Let $f:\g\to\h$ be a morphism of finite-dimensional real Lie triple systems. If $H_f^2=0$, then $\mathfrak U(f)$ is graded stable under simultaneous perturbations of the two enveloping Lie brackets. For smooth families of such brackets through the original ones, a nearby graded homomorphism can be chosen smoothly through $\mathfrak U(f)$.
\end{theorem}
\begin{proof}
By the recognition theorem for universal envelopes \cite{Smirnov}, a graded Lie algebra generated by its odd component is the universal envelope of that odd Lie triple system when its second graded Chevalley--Eilenberg cohomology with every trivial even coefficient module vanishes. For the source universal envelope, the corresponding conditions are
\begin{equation}\label{eq:envelope-recognition}
 [\mathfrak U(\g)_1,\mathfrak U(\g)_1]=\mathfrak U(\g)_0,
 \qquad H_{\mathrm{CE},\mathrm{gr}}^2
       (\mathfrak U(\g),\R_{\mathrm{even}})=0.
\end{equation}
The coefficient module here is trivial and concentrated in even degree. Both conditions hold for $\mathfrak U(\g)$. Put $L=\mathfrak U(\g)$ and keep the decomposition $L=L_0\oplus L_1$ fixed. For every trivial even coefficient module $M$, finite dimensionality of $L$ gives
\[
 C_{\mathrm{CE},\mathrm{gr}}^n(L,M)
 \cong C_{\mathrm{CE},\mathrm{gr}}^n(L,\R_{\mathrm{even}})\otimes M.
\]
Under this identification the differential is $d^n\otimes\id_M$. Tensoring over $\R$ preserves exactness, so vanishing for $\R_{\mathrm{even}}$ gives vanishing for every such $M$.

Both conditions persist under sufficiently small graded Lie bracket perturbations on the fixed vector space $L$. The first is surjectivity of the bracket map from $\wedge^2L_1$ to $L_0$, which is open. Hence every sufficiently close graded Lie bracket still satisfies $[L_1,L_1]'=L_0$, so its odd component generates the whole Lie algebra.

For the second, consider the fixed finite-dimensional spaces
\[
 C_{\mathrm{CE},\mathrm{gr}}^n(L,\R_{\mathrm{even}})
 =\Hom\bigl((\wedge^nL)_0,\R\bigr),\qquad n=1,2,3.
\]
Here $(\wedge^nL)_0$ is the even component of the ordinary exterior power. Write $d^1,d^2$ for the original differentials and $d^{\prime 1},d^{\prime 2}$ for those obtained from the perturbed bracket. In fixed bases, their matrix entries depend linearly on the bracket coefficients. Nonzero minors realizing the original ranks remain nonzero in a common neighborhood. Thus neither rank can decrease there. Exactness at the original bracket gives
\[
 \operatorname{rank}d^1+\operatorname{rank}d^2
 =\dim C_{\mathrm{CE},\mathrm{gr}}^2(L,\R_{\mathrm{even}}).
\]
For every nearby bracket satisfying the Jacobi identity, $d^{\prime 2}d^{\prime 1}=0$, and therefore
\begin{align*}
 \dim C_{\mathrm{CE},\mathrm{gr}}^2(L,\R_{\mathrm{even}})
 &=\operatorname{rank}d^1+\operatorname{rank}d^2\\
 &\leq\operatorname{rank}d^{\prime 1}+\operatorname{rank}d^{\prime 2}
 \leq\dim C_{\mathrm{CE},\mathrm{gr}}^2(L,\R_{\mathrm{even}}).
\end{align*}
Equality follows, so $\im d^{\prime 1}=\ker d^{\prime 2}$. The same tensor product argument gives the required vanishing for every trivial even coefficient module on this common neighborhood. Thus all hypotheses of the recognition theorem hold for the nearby graded source brackets.

More explicitly, let $L'$ denote $L$ with such a bracket, and let $T'=L_1$ have the triple product $[x,y,z]'=[[x,y]',z]'$. The identity on $L_1$ induces a canonical graded homomorphism
\[
 p':\mathfrak U(T')\longrightarrow L',\qquad
 p'(x)=x,\quad p'([x,y])=[x,y]'.
\]
It is surjective because $[L_1,L_1]'=L_0$. Its kernel lies in the even component and commutes with all odd generators: the bracket of a kernel element with an odd generator is odd and is mapped to zero, while $p'$ is the identity on the odd component. Hence $\ker p'$ is central. Over $\R$, this surjection admits a grading-preserving linear section. The vanishing of $H_{\mathrm{CE},\mathrm{gr}}^2(L',\ker p')$ then gives a graded Lie splitting $s'$ with $p's'=\id_{L'}$, by the central extension characterization in \cite{Smirnov}. Since $s'$ is the identity on the odd component, its image contains all generators of $\mathfrak U(T')$. Thus $s'$ is surjective, and $p'$ is an isomorphism restricting to the identity on $T'$. This identifies the perturbed source with the universal envelope of its own odd triple system and retains the universal property relative to the fixed odd inclusion.

Now perturb the source and target graded brackets on the fixed spaces $\mathfrak U(\g)$ and $\mathfrak U(\h)$. Their iterated brackets give nearby triple products on $\g$ and $\h$. Theorem~\ref{thm:stability} supplies a nearby odd map $q:\g\to\h$ satisfying
\begin{equation*}
 q\bigl([[x,y]'_{\mathfrak U(\g)},z]'_{\mathfrak U(\g)}\bigr) =[[qx,qy]'_{\mathfrak U(\h)},qz]'_{\mathfrak U(\h)}.
\end{equation*}
The universal property first extends $q$ to a graded Lie homomorphism $\widetilde q:\mathfrak U(T')\to(\mathfrak U(\h),[\ ,\ ]')$. The required map on the fixed perturbed source is $\widehat q=\widetilde q\circ(p')^{-1}$. It restricts to $q$ on the odd component and is unique because the odd component generates the source. At the original brackets, the extension is $\mathfrak U(f)$.

To see its closeness and smooth dependence, choose $x_i,y_i\in\g$ whose original brackets form an even basis. Their perturbed brackets $[x_i,y_i]'_{\mathfrak U(\g)}$ remain a basis, and the extension must take their values to
\[
 [qx_i,qy_i]'_{\mathfrak U(\h)}.
\]
More explicitly, for $a\in L_0$ write $a=\sum_i c_i'[x_i,y_i]'_{\mathfrak U(\g)}$. The coordinates $c_i'$ depend linearly on $a$ and smoothly on the source bracket, by inversion of the matrix of this moving basis. For $x\in\g$, the extension is therefore
\[
 \widehat q(a+x)=\sum_i c_i'[qx_i,qy_i]'_{\mathfrak U(\h)}+qx.
\]
The basis change, the target brackets and the odd map supplied by Theorem~\ref{thm:stability} depend smoothly on the parameters. Hence the whole graded extension depends smoothly on the perturbed brackets and equals $\mathfrak U(f)$ at the base point. This proves graded stability and the smooth-family assertion. Here the recognition step uses graded Chevalley--Eilenberg cohomology with trivial even coefficients, while $H_f^2$ controls the morphism stability problem.
\end{proof}

Theorems~\ref{thm:graded-local-comparison} and~\ref{thm:graded-stability} compare the local Lie triple problem with grading-preserving Lie deformations. The corresponding inner equivalences are generated by even Lie inner derivations.

For prescribed ungraded Lie algebras, the comparison takes a different form.

\begin{theorem}\label{thm:ordinary-local-comparison}
Let $f:\mathfrak a\to\mathfrak b$ be a homomorphism of finite-dimensional real Lie algebras satisfying Eq.~\eqref{eq:Lie-comparison-hyp}. Then the spaces of Lie algebra homomorphisms and associated Lie triple morphisms coincide in a neighborhood of $f$. This coincidence remains valid, uniformly near $f$, when both binary Lie brackets are perturbed. Consequently, stability of $f$ in the Lie category is equivalent to stability of the associated triple morphism along the family of triple products induced by these binary perturbations.

If $\mathfrak b$ is also perfect, local inner rigidity in the Lie and Lie triple categories is equivalent. Moreover, $H_f^2=0$ implies $H_{\mathrm{CE}}^2(f)=0$ and Lie stability. Under perfectness of $\mathfrak b$, the condition $H_f^1=0$ is equivalent to $H_{\mathrm{CE}}^1(f)=0$ and implies local inner rigidity in both categories.
\end{theorem}
\begin{proof}
Lemma~\ref{lem:T-injective} makes the linear map
\[
 b\longmapsto\bigl((x,y,z)\longmapsto
       [b(x,y),fz]_{\mathfrak b}+b([x,y]_{\mathfrak a},z)\bigr)
\]
injective. Injectivity is open for linear maps between fixed finite-dimensional spaces. Thus, for all nearby $q$ and nearby binary Lie brackets, the map
\[
 b\longmapsto\bigl((x,y,z)\longmapsto
       [b(x,y),qz]'_{\mathfrak b}+b([x,y]'_{\mathfrak a},z)\bigr)
\]
remains injective in a common neighborhood. Expanding the identity from Proposition~\ref{prop:binary-defect} with the perturbed brackets gives
\begin{equation*}
 [[qx,qy]'_{\mathfrak b},qz]'_{\mathfrak b} -q[[x,y]'_{\mathfrak a},z]'_{\mathfrak a} =\bigl[[qx,qy]'_{\mathfrak b}-q[x,y]'_{\mathfrak a}, qz\bigr]'_{\mathfrak b} +[q[x,y]'_{\mathfrak a},qz]'_{\mathfrak b} -q[[x,y]'_{\mathfrak a},z]'_{\mathfrak a}.
\end{equation*}
By the displayed injectivity, the left-hand side vanishes for all $x,y,z$ if and only if
\[
 [qx,qy]'_{\mathfrak b}-q[x,y]'_{\mathfrak a}=0
 \qquad(x,y\in\mathfrak a).
\]
This proves the uniform local equality of the morphism spaces and the equivalence of stability along the induced family of triple products.

The source is perfect by hypothesis. If the target is also perfect, every Lie inner action on either space is a sum of actions $z\mapsto[[u,v],z]$, which are exactly the Lie triple inner actions. Thus their inner orbits agree, proving equivalence of local inner rigidity. Proposition~\ref{prop:ordinary-cohomology} identifies the first cohomology groups under these hypotheses and gives an injection from second Lie cohomology to second Lie triple cohomology. Hence $H_f^2=0$ implies $H_{\mathrm{CE}}^2(f)=0$. Lie stability also follows directly by applying Theorem~\ref{thm:stability} to the induced triple products and then using the local equality just proved. Finally, the vanishing of the first cohomology implies local inner rigidity by Theorem~\ref{thm:local-rigid}.
\end{proof}

This stability equivalence concerns perturbations of triple products induced by binary Lie brackets. The injective map in Eq.~\eqref{eq:second-comparison} carries Lie obstruction classes to Lie triple obstruction classes.

\begin{example}\label{ex:Heisenberg-comparison}
Let $\mathfrak n$ be the real Heisenberg Lie algebra with basis $p,q,z$ and bracket $[p,q]=z$. Its associated triple product is zero. Hence the maps
\[
 f_t(p)=p,\qquad f_t(q)=q,\qquad f_t(z)=(1+t)z
\]
form a smooth Lie triple deformation of $\id_{\mathfrak n}$. For $t\neq0$, they are not Lie algebra homomorphisms, since $[f_t(p),f_t(q)]=z\neq f_t([p,q])$. Thus the local Lie morphism space is a proper subset of the local Lie triple morphism space in this example. All Lie triple inner derivations are zero, whereas the Lie inner derivations $\ad_p$ and $\ad_q$ are nonzero, giving distinct boundary spaces in Eq.~\eqref{eq:ordinary-boundaries}.
\end{example}

\section{Triple subsystems and the graph comparison}

In this section, we study triple subsystems through their normal representations and graph defects. We establish local rigidity, smoothness and stability, and compare their cohomology with morphism cohomology. Finally, we relate triple graph subsystems to graded Lie graph subalgebras and compare the normal differentials for Lie subalgebras and their associated triple subsystems.

Let $E$ be a Lie triple system and let $\s\subseteq E$ be a triple subsystem. We equip the quotient vector space $E/\s$ with the normal representation defined below. Write $\overline u=u+\s$.

\begin{proposition}\label{prop:normal-rep}
The maps
\begin{equation}\label{eq:normal-rep}
 \theta_N(x,y)\overline u=\overline{[u,x,y]},\qquad
 D_N(x,y)\overline u=\overline{[x,y,u]}
 \quad(x,y\in\s)
\end{equation}
define a representation of $\s$ on $E/\s$.
\end{proposition}
\begin{proof}
For $x,y,v\in\s$ and $u\in E$, closure of $\s$ under the triple product gives
\begin{align*}
 [u+v,x,y]+\s&=[u,x,y]+[v,x,y]+\s=[u,x,y]+\s,\\
 [x,y,u+v]+\s&=[x,y,u]+[x,y,v]+\s=[x,y,u]+\s.
\end{align*}
Thus the two quotient operators are well defined. Restrict the adjoint representation of $E$ to $\s$. Every operator $u\mapsto[u,x,y]$ with $x,y\in\s$ preserves $\s$, so the adjoint representation identities descend to $E/\s$. Moreover, the cyclic identity gives
\[
 [u,y,x]-[u,x,y]+\s=[x,y,u]+\s,
\]
which is the required identity between $\theta_N$ and $D_N$.
\end{proof}

For $n\geq1$, put
\[
 C_N^n(\s,E)=\C^{n-1}(\s,E/\s).
\]
Let $\partial^n$ be the Yamaguti differential with coefficients \eqref{eq:normal-rep}, with the same indexing as \eqref{eq:full-differential}. Denote the resulting positive-degree groups by $H_N^n(\s,E)$, where
\[
 H_N^1(\s,E)=Z_N^1(\s,E)=\ker\partial^1.
\]
In particular,
\begin{align}
 (\partial^1a)(x,y,z)
 &=\theta_N(y,z)a(x)-\theta_N(x,z)a(y)
      +D_N(x,y)a(z)-a([x,y,z]). \label{eq:normal-d1}
\end{align}
For comparison with inner orbits, augment this complex by
\begin{equation}\label{eq:normal-augmentation}
 \widetilde C_N^0(\s,E)=\wedge^2E,
 \qquad (\partial^0X)(x)=\overline{\ad_Xx}.
\end{equation}
A derivation calculation as in Theorem~\ref{thm:complex} gives $\partial^1\partial^0=0$. Write $\widetilde H_N^n$ for this augmented cohomology. Thus
\[
 \widetilde H_N^1=Z_N^1/\im\partial^0,
 \qquad \widetilde H_N^n=H_N^n\quad(n\geq2).
\]

Suppose now that $E$ is finite-dimensional and real, and let $k=\dim\s$. A deformation of $\s$ is a smooth curve of $k$-dimensional triple subsystems in $\Gr_k(E)$ through $\s$. Choose a vector-space complement $E=\s\oplus\m$ and identify $\m$ with $E/\s$. Every nearby $k$-plane is uniquely the graph
\[
 \s_a=\{x+a(x):x\in\s\},\qquad a\in\Hom(\s,\m).
\]
For a trilinear tensor $\nu$ with Lie triple symmetries, define
\begin{equation}\label{eq:subsystem-defect}
 S_\nu(a)(x,y,z) =\pr_\m\nu(x+ax,y+ay,z+az) -a\pr_\s\nu(x+ax,y+ay,z+az).
\end{equation}
This is a cochain in $\C^1(\s,\m)$. It vanishes exactly when $\s_a$ is closed under $\nu$.

\begin{proposition}\label{prop:subsystem-inf}
The derivative of a smooth subsystem deformation belongs to $Z_N^1(\s,E)$. More precisely, under the graph chart,
\begin{equation}\label{eq:subsystem-linearization}
 (dS_\nu)_0=\partial^1
\end{equation}
when $\nu$ is the original law of $E$ and $\s$ is its subsystem.
\end{proposition}
\begin{proof}
For $a\in\Hom(\s,\m)$ and $x,y,z\in\s$, differentiation of \eqref{eq:subsystem-defect} at zero gives
\begin{align*}
 (dS_\nu)_0(a)(x,y,z)
 &=\pr_\m\bigl([ax,y,z]+[x,ay,z]+[x,y,az]\bigr)
       -a([x,y,z]).
\end{align*}
Under the identification of $\m$ with $E/\s$, this is the expression in \eqref{eq:normal-d1}. If the graphs of $a_t$ are subsystems and $a_0=0$, then
\begin{equation*}
 \pr_\m[x+a_tx,y+a_ty,z+a_tz] -a_t\pr_\s[x+a_tx,y+a_ty,z+a_tz]=0.
\end{equation*}
Differentiating at $t=0$ yields
\[
 \pr_\m\bigl([\dot a_0x,y,z]+[x,\dot a_0y,z]
       +[x,y,\dot a_0z]\bigr)-\dot a_0([x,y,z])=0.
\]
Hence $\partial^1\dot a_0=0$. The tangent space of the Grassmannian is canonically $\Hom(\s,E/\s)$, so the conclusion is independent of the complement.
\end{proof}

We construct the compatibility map for \eqref{eq:subsystem-defect} explicitly. Let
\[
 \mathsf P_\s:\C^\bullet(E,E)\longrightarrow\C^\bullet(\s,\m)
\]
be restriction to $\s$ inputs followed by projection to $\m$, and extend $\s$-to-$\m$ cochains by zero whenever an input lies in $\m$. This subspace of cochains is abelian, and $\ker\mathsf P_\s$ is a graded Lie subalgebra, by the same argument as in Proposition~\ref{prop:derived}. Extend $a$ as an endomorphism $\widehat a$ of $E$ with $\widehat a|_\m=0$, and set
\begin{align}
 \nu_a&=(\id-\widehat a)\nu(\id+\widehat a,\id+\widehat a,\id+\widehat a),
       \label{eq:subsystem-conjugate}\\
 \mathcal B_{\nu,a}(b)&=\mathsf P_\s[\nu_a,\widehat b]_{\B},
       \qquad b\in\C^1(\s,\m). \label{eq:subsystem-B}
\end{align}

\begin{lemma}\label{lem:subsystem-bianchi}
If $\nu$ is a Lie triple system law on $E$, then
\begin{equation}\label{eq:subsystem-bianchi}
 \mathcal B_{\nu,a}\bigl(S_\nu(a)\bigr)=0
\end{equation}
for every $a$. At the original law and $a=0$, one has $\mathcal B_{\nu,0}=-\partial^2$.
\end{lemma}
\begin{proof}
Conjugation by $\id+\widehat a$ preserves the Lie triple identities. Since $\ker\mathsf P_\s$ is a graded Lie subalgebra and the subspace of extended $\s$-to-$\m$ cochains is abelian, projecting the graded self-bracket gives
\begin{align*}
 0={}&\mathsf P_\s\Bigl[ (\id-\widehat a)\nu(\id+\widehat a,\id+\widehat a,\id+\widehat a),\qquad (\id-\widehat a)\nu(\id+\widehat a,\id+\widehat a,\id+\widehat a) \Bigr]_{\B} \\
 ={}&2\mathsf P_\s\Bigl[ (\id-\widehat a)\nu(\id+\widehat a,\id+\widehat a,\id+\widehat a),\qquad \widehat{\mathsf P_\s\bigl((\id-\widehat a) \nu(\id+\widehat a,\id+\widehat a,\id+\widehat a)\bigr)} \Bigr]_{\B}.
\end{align*}
On inputs $x,y,z\in\s$, the projected trilinear map in the second argument is
\[
 \pr_\m\nu(x+ax,y+ay,z+az)
       -a\pr_\s\nu(x+ax,y+ay,z+az).
\]
Thus the last equality is precisely \eqref{eq:subsystem-bianchi}, by \eqref{eq:subsystem-B}.

At the original law and $a=0$, direct insertion of $b\in\C^1(\s,\m)$ gives
\begin{align*}
 &(\mathsf P_\s[\nu,\widehat b]_{\B})(x,y,u,v,w) =\pr_\m[b(x,y,u),v,w]+\pr_\m[u,b(x,y,v),w] \\
 &\quad+\pr_\m[u,v,b(x,y,w)]-\pr_\m[x,y,b(u,v,w)] +b([x,y,u],v,w)+b(u,[x,y,v],w) \\
 &\quad+b(u,v,[x,y,w])-b(x,y,[u,v,w]) =-(\partial^2b)(x,y,u,v,w).
\end{align*}
This proves the asserted value at $a=0$. More generally, the same insertion formula in degree $n-1$, multiplied by $(-1)^{n-1}$, is the Yamaguti differential for the quotient representation \eqref{eq:normal-rep}.
\end{proof}

\begin{remark}\label{rem:subsystem-derived}
If $\s$ is a subsystem for $\nu$, then $\mathsf P_\s\nu=0$ and the same data give a derived $L_\infty[1]$-algebra on $\C^\bullet(\s,\m)$. Its brackets are
\[
 m_j(b_1,\ldots,b_j)
   =\mathsf P_\s[\cdots[[\nu,\widehat b_1]_{\B},\widehat b_2]_{\B},\ldots,\widehat b_j]_{\B}.
\]
For a degree-zero map $a$, its Maurer--Cartan expression is precisely $S_\nu(a)$. In general the expression has degree four in $a$; its fourth-degree term is
\[
 -a\pr_\s\nu(ax,ay,az).
\]
For morphism graphs, the ambient law is a direct product and the complementary target summand is an ideal. Hence $\pr_\s\nu(ax,ay,az)=0$, and the graph equation is cubic. Formula \eqref{eq:subsystem-conjugate} applies to every vector-space complement $\m$.
\end{remark}

\begin{definition}
A $k$-dimensional subsystem $\s$ of $(E,\nu)$ is stable if, for every neighborhood $U$ of $\s$ in $\Gr_k(E)$, every sufficiently close Lie triple system law $\nu'$ on $E$ has a $k$-dimensional subsystem in $U$.
\end{definition}

\begin{theorem}\label{thm:subsystem-stable}
If $H_N^2(\s,E)=0$, the space of $k$-dimensional triple subsystems of $E$ is a smooth submanifold near $\s$, of dimension $\dim Z_N^1(\s,E)$. Every normal first cocycle is tangent to a smooth subsystem deformation. Moreover, $\s$ is stable, and for a smooth family of nearby ambient laws one can choose a smooth family of nearby subsystems.
\end{theorem}
\begin{proof}
In the Grassmannian graph chart, apply Lemma~\ref{lem:exact-zero} to the polynomial map from $\Hom(\s,\m)$ to $\C^1(\s,\m)$ given by
\begin{equation*}
 a\longmapsto\bigl((x,y,z)\longmapsto \pr_\m\nu(x+ax,y+ay,z+az) -a\pr_\s\nu(x+ax,y+ay,z+az)\bigr),
\end{equation*}
and to the family of linear maps from $\C^1(\s,\m)$ to $\C^2(\s,\m)$ given by
\[
 b\longmapsto\mathsf P_\s\bigl[
 (\id-\widehat a)\nu(\id+\widehat a,\id+\widehat a,\id+\widehat a),
 \widehat b\bigr]_{\B}.
\]
Lemma~\ref{lem:subsystem-bianchi} supplies their compatibility. The derivative of the polynomial map at zero sends $a$ to
\[
 (x,y,z)\longmapsto
 \pr_\m\bigl([ax,y,z]+[x,ay,z]+[x,y,az]\bigr)-a([x,y,z]),
\]
which is $\partial^1a$. The family of linear maps at zero sends $b$ to $-\partial^2b$. Thus the required exactness is $\im\partial^1=\ker\partial^2$, which is the hypothesis $H_N^2(\s,E)=0$. The lemma gives local smoothness, tangent space $\ker\partial^1$ and dimension $\dim Z_N^1(\s,E)$. Every tangent vector is realized by a smooth curve in this local submanifold.

For stability, allow $\nu'$ to range over all trilinear tensors with the Lie triple symmetries. The two displayed maps, with $\nu$ replaced by $\nu'$, depend polynomially on $(\nu',a)$. When $\nu'$ satisfies the fundamental identity, Lemma~\ref{lem:subsystem-bianchi} supplies compatibility. The parameter part of Lemma~\ref{lem:exact-zero} gives a smooth choice $a(\nu')$, with $a(\nu)=0$, such that
\begin{equation*}
 \pr_\m\nu'\bigl(x+a(\nu')x,y+a(\nu')y,z+a(\nu')z\bigr) =a(\nu')\pr_\s\nu' \bigl(x+a(\nu')x,y+a(\nu')y,z+a(\nu')z\bigr).
\end{equation*}
Its graph is the required nearby subsystem. Continuity gives stability, and composition with a smooth curve of ambient laws gives a smooth curve of subsystems.
\end{proof}

\begin{theorem}\label{thm:subsystem-rigid}
If $\widetilde H_N^1(\s,E)=0$, every sufficiently close $k$-dimensional triple subsystem is the image of $\s$ under an inner automorphism of $E$. The automorphism can be chosen smoothly in a neighborhood of $\s$ along the subsystem locus.
\end{theorem}
\begin{proof}
For $X\in\wedge^2E$, the derivative at zero of $t\mapsto\exp(t\ad_X)\s$ is
\[
 x\longmapsto\ad_Xx+\s.
\]
Thus the tangent space to the inner orbit is $\im\partial^0$. The hypothesis $\widetilde H_N^1(\s,E)=0$ identifies this space with $\ker\partial^1$.

Choose a complement to $\im\partial^0$ in $\Hom(\s,\m)$ and a complement to the kernel of the orbit derivative in $\InnDer(E)$. The action of exponentials on graph planes and the inverse function theorem give local coordinates consisting of an orbit coordinate and a transverse graph coordinate. Inner automorphisms preserve the subsystem condition. On the transverse complement, the linear map
\[
 a\longmapsto\bigl((x,y,z)\longmapsto
 \pr_\m([ax,y,z]+[x,ay,z]+[x,y,az])-a([x,y,z])\bigr)
\]
is injective. The finite-dimensional Taylor estimate used in Theorem~\ref{thm:local-rigid} therefore yields constants $c,C>0$ such that, for sufficiently small transverse $a$,
\begin{equation*}
 \bigl\|\pr_\m\nu(\cdot+a\,\cdot,\cdot+a\,\cdot,\cdot+a\,\cdot) -a\pr_\s\nu(\cdot+a\,\cdot,\cdot+a\,\cdot,\cdot+a\,\cdot) \bigr\|\geq c\|a\|-C\|a\|^2.
\end{equation*}
Hence a sufficiently small transverse graph is a subsystem only when $a=0$. Every nearby subsystem is consequently the image of $\s$ under an inner automorphism. The inverse coordinate map is smooth, and its orbit coordinate gives the claimed smooth choice of inner automorphism.
\end{proof}

\begin{proposition}\label{prop:graph-germ}
Let $f:(\g,\pi)\to(\h,\mu)$ be a morphism and set
\[
 G_f=\{(x,fx):x\in\g\}\subseteq \g\oplus\h.
\]
Then $G_f$ is a triple subsystem of the direct product. In finite dimension, the graph map identifies the local space of morphisms near $f$ with the local space of $\dim\g$-dimensional triple subsystems near $G_f$.
\end{proposition}
\begin{proof}
The direct product of $(x,fx),(y,fy),(z,fz)$ is
\[
 (\pi(x,y,z),\mu(fx,fy,fz))
 =(\pi(x,y,z),f\pi(x,y,z)),
\]
which belongs to the graph of $f$. Conversely, membership in this graph requires $\mu(fx,fy,fz)=f\pi(x,y,z)$. The graph of $f$ is transverse to $0\oplus\h$. In finite dimension this transversality is open, so every sufficiently close plane of dimension $\dim\g$ is the graph of a unique linear map $q:\g\to\h$. Closure of that graph under the triple product is exactly
\[
 \mu(qx,qy,qz)=q\pi(x,y,z).
\]
The Grassmannian graph chart and its inverse are smooth, proving the local identification.
\end{proof}

There are canonical vector-space identifications
\begin{equation}\label{eq:graph-identification}
 j_f:\g\longrightarrow G_f,\quad j_f(x)=(x,fx),\qquad
 Q_f:(\g\oplus\h)/G_f\longrightarrow\h,\quad
 Q_f(\overline{(x,u)})=u-fx.
\end{equation}

\begin{theorem}\label{thm:graph-cohomology}
For every $n\geq1$, define
\begin{align}
 &(\Xi_n a)(j_fx_1,\ldots,j_fx_{2n-1})
     =\overline{(0,a(x_1,\ldots,x_{2n-1}))}. \label{eq:Xi}
\end{align}
Then $\Xi_n$ identifies the positive-degree morphism complex with the normal complex of $G_f$. It also identifies their first boundary spaces when the normal complex is augmented by $\wedge^2(\g\oplus\h)$. Consequently,
\begin{equation}\label{eq:graph-aug-iso}
 H_f^n\cong\widetilde H_N^n(G_f,\g\oplus\h)
 \qquad(n\geq1).
\end{equation}
For the unaugmented normal complex, the comparison is
\begin{align}
 H_f^n&\cong H_N^n(G_f,\g\oplus\h)\quad(n\geq2), \label{eq:graph-high}\\
 0\longrightarrow B_f^1&\xrightarrow{\ \Xi_1\ }
 H_N^1(G_f,\g\oplus\h)\longrightarrow H_f^1\longrightarrow0.
 \label{eq:graph-exact}
\end{align}
\end{theorem}
\begin{proof}
The identification $(x,fx)\mapsto x$ of the graph with $\g$ and the quotient identification $(z,u)+G_f\mapsto u-fz$ show that \eqref{eq:Xi} is an isomorphism in each positive degree. These identifications preserve the cochain conditions in the last three variables. For $x,y,z\in\g$ and $u\in\h$, the normal coefficient action is
\begin{align*}
 \theta_N((x,fx),(y,fy))\overline{(z,u)} &=\overline{(\pi(z,x,y),\mu(u,fx,fy))} =\overline{(0,\mu(u,fx,fy)-f\pi(z,x,y))} \\
 &=\overline{(0,\mu(u-fz,fx,fy))}.
\end{align*}
Thus the quotient identification carries this normal action to the coefficient action in \eqref{eq:pullback}. The map $x\mapsto(x,fx)$ also preserves the source triple product. Substituting these two facts into each term of \eqref{eq:full-differential} gives
\begin{align*}
 &(\partial^n\Xi_na)((x_1,fx_1),\ldots,(x_{2n+1},fx_{2n+1})) =\overline{(0,(\delta_f^na)(x_1,\ldots,x_{2n+1}))} \\
 &\qquad=(\Xi_{n+1}\delta_f^na) ((x_1,fx_1),\ldots,(x_{2n+1},fx_{2n+1})).
\end{align*}
Hence the positive-degree complexes are identified.

For the incoming maps in degree one, use the canonical decomposition
\[
 \wedge^2(\g\oplus\h)
 =\wedge^2\g\oplus(\g\otimes\h)\oplus\wedge^2\h.
\]
Mixed pairs act trivially in the direct product. For $X\in\wedge^2\g$ and $U\in\wedge^2\h$, the resulting normal orbit derivative sends $(x,fx)$ to
\begin{equation*}
 \overline{(\ad_Xx,\ad_Ufx)} =\overline{(0,\ad_Ufx-f\ad_Xx)} =\overline{(0,(\delta_f^0(X,U))(x))}.
\end{equation*}
Every ambient inner derivative has this form, and every pair $(X,U)$ is represented by an ambient two-vector. Thus $\Xi_1$ identifies $\im\delta_f^0$ with $\im\partial^0$, proving \eqref{eq:graph-aug-iso}. The first unaugmented normal cohomology is the first cocycle space itself. Quotienting by the displayed boundary space gives \eqref{eq:graph-exact}. The groups in degrees at least two are unchanged, proving \eqref{eq:graph-high}.
\end{proof}

\begin{example}
For the identity of the two-dimensional system in Example~\ref{ex:rigid}, $H_f^1=0$, while
\[
 H_N^1(G_f,T\oplus T)\cong Z_f^1=\R J.
\]
Thus the first unaugmented normal group has dimension one. Its quotient by the one-dimensional ambient inner orbit tangent space is zero, in agreement with \eqref{eq:graph-aug-iso}.
\end{example}

Theorem~\ref{thm:graph-cohomology} and Proposition~\ref{prop:graph-germ} show that the morphism and graph-subsystem rigidity criteria agree after augmentation. Their smoothness and stability criteria agree through the degree-two cohomology isomorphism. Under a perturbation of the two direct product factors, a nearby graph subsystem remains transverse to the target summand, so it again determines a morphism. For a general ambient perturbation of $\g\oplus\h$, Theorem~\ref{thm:subsystem-stable} gives a nearby triple subsystem under the same degree-two vanishing condition.

We finish by relating the graph construction to Lie subalgebras and Lie algebra homomorphisms. For a triple subsystem $\s\subseteq E$, its inclusion induces a Lie homomorphism $\mathfrak U(\s)\to\mathfrak U(E)$. Its image is the graded Lie subalgebra generated by $\s$, namely
\begin{equation}\label{eq:subsystem-lie-generated}
 [\s,\s]_{\mathfrak U(E)}\oplus\s.
\end{equation}
The kernel is contained in the even component and is central: it acts trivially on the odd generators, and these generate $\mathfrak U(\s)$. Thus the generated subalgebra in Eq.~\eqref{eq:subsystem-lie-generated} is the quotient of $\mathfrak U(\s)$ by this central kernel. For morphism graphs, the generated subalgebra has the following explicit description.

\begin{proposition}\label{prop:enveloping-graph}
Let $\Phi=\mathfrak U(f)$, and let
\[
 \Gamma_\Phi=\{(X,\Phi X):X\in\mathfrak U(\g)\}
       \subseteq\mathfrak U(\g)\oplus\mathfrak U(\h).
\]
Then $\Gamma_\Phi$ is the graded Lie subalgebra generated by its odd component, and
\begin{equation}\label{eq:odd-graph}
 (\Gamma_\Phi)_1=G_f,\qquad
 (\Gamma_\Phi)_0
 =\operatorname{span}\{([x,y],[fx,fy]):x,y\in\g\}.
\end{equation}
In finite dimension, the map $G_q\mapsto\Gamma_{\mathfrak U(q)}$ identifies the local space of triple graph subsystems with the local space of graded Lie graph subalgebras. In particular, the graded Lie subalgebra germs considered here have the fixed even and odd dimensions of $\Gamma_\Phi$ and are transverse to the target summand.
\end{proposition}
\begin{proof}
For $X,Y\in\mathfrak U(\g)$, preservation of the binary bracket gives
\begin{equation*}
 [(X,\mathfrak U(f)X),(Y,\mathfrak U(f)Y)] =([X,Y],[\mathfrak U(f)X,\mathfrak U(f)Y]) =([X,Y],\mathfrak U(f)[X,Y]).
\end{equation*}
Thus the graph is a Lie subalgebra, and it is graded because $\mathfrak U(f)$ preserves parity. Its odd component consists of $(x,fx)$ with $x\in\g$. Bracketing two such elements gives
\[
 [(x,fx),(y,fy)]=([x,y],[fx,fy]).
\]
These brackets span the even component, since the odd generators generate $\mathfrak U(\g)$. This proves \eqref{eq:odd-graph} and the generation assertion.

In finite dimension, a sufficiently close graded plane of the same graded dimensions, transverse to the target summand, is the graph of a unique parity-preserving linear map. Closure of the graph under the bracket is exactly preservation of that bracket. Proposition~\ref{prop:universal-lie} identifies these graded Lie homomorphisms with the universal extensions of Lie triple morphisms. Proposition~\ref{prop:graph-germ} identifies their odd graphs with the nearby triple graph subsystems. The polynomial extension and linear restriction from Theorem~\ref{thm:graded-local-comparison}, together with the smooth Grassmannian graph charts, give the stated local identification.
\end{proof}

Consequently, formal and local inner rigidity of a triple graph, with the augmentation in Eq.~\eqref{eq:normal-augmentation}, corresponds to inner rigidity of its graded Lie graph under even ambient inner automorphisms. Projection onto the source identifies the even graph with $\mathfrak U(\g)_0$, fixing its dimension. Stability of $\mathfrak U(f)$ from Theorem~\ref{thm:graded-stability} also gives persistence of its graph under small graded direct product Lie bracket perturbations.

There is a further direct comparison when the triple systems come from prescribed Lie algebras. Let $\mathfrak e$ be a Lie algebra and $\mathfrak s\subseteq\mathfrak e$ a Lie subalgebra. The normal Lie representation is
\[
 \rho_N(x)\overline u=\overline{[x,u]},
 \qquad x\in\mathfrak s,\quad \overline u\in\mathfrak e/\mathfrak s.
\]
Its Chevalley--Eilenberg differential in degree one is
\begin{equation}\label{eq:Lie-normal-d1}
 (d_Na)(x,y)=\rho_N(x)a(y)-\rho_N(y)a(x)-a([x,y]).
\end{equation}
For an alternating map $b:\wedge^2\mathfrak s\to\mathfrak e/\mathfrak s$, set
\[
 (\mathcal T_Nb)(x,y,z)=b([x,y],z)-\rho_N(z)b(x,y).
\]

\begin{proposition}\label{prop:Lie-normal-comparison}
For the associated Lie triple subsystem, the normal representation satisfies
\begin{equation}\label{eq:Lie-triple-normal-action}
 \theta_N(x,y)=\rho_N(y)\rho_N(x),\qquad
 D_N(x,y)=\rho_N([x,y]).
\end{equation}
Its first coboundary operator is related to the Lie normal differential by
\begin{equation}\label{eq:normal-differential-comparison}
 \partial^1a=\mathcal T_N(d_Na).
\end{equation}
Thus every infinitesimal Lie subalgebra deformation is an infinitesimal triple subsystem deformation. For a Lie homomorphism $f:\mathfrak a\to\mathfrak b$, the ordinary Lie normal complex of its graph is canonically the Chevalley--Eilenberg morphism complex of $f$, and Eq.~\eqref{eq:normal-differential-comparison} becomes Eq.~\eqref{eq:differential-comparison} under the graph identifications.
\end{proposition}
\begin{proof}
For $x,y\in\mathfrak s$ and $u\in\mathfrak e$, the normal Lie action gives
\begin{equation*}
 \rho_N(y)\rho_N(x)(u+\mathfrak s) =[y,[x,u]]+\mathfrak s=[[u,x],y]+\mathfrak s,\qquad \rho_N([x,y])(u+\mathfrak s) =[[x,y],u]+\mathfrak s.
\end{equation*}
These are exactly the two normal Lie triple actions in \eqref{eq:normal-rep}. Expanding the right-hand side of \eqref{eq:normal-differential-comparison} gives
\begin{align*}
 &\rho_N([x,y])a(z)-\rho_N(z)a([x,y])-a([[x,y],z]) -\rho_N(z)\bigl(\rho_N(x)a(y)-\rho_N(y)a(x)-a([x,y])\bigr) \\
 &=\rho_N(z)\rho_N(y)a(x)-\rho_N(z)\rho_N(x)a(y) +\rho_N([x,y])a(z)-a([[x,y],z]).
\end{align*}
The two action identities just proved identify this expression with
\[
 \theta_N(y,z)a(x)-\theta_N(x,z)a(y)
             +D_N(x,y)a(z)-a([[x,y],z]),
\]
which is the first normal Lie triple differential. In particular, a vanishing first Lie normal differential gives a vanishing first Lie triple normal differential.

For the graph of $f$, the normal Lie action is explicitly
\begin{equation*}
 \rho_N((x,fx))\overline{(y,u)} =\overline{([x,y]_{\mathfrak a},[fx,u]_{\mathfrak b})} =\overline{(0,[fx,u]_{\mathfrak b}-f[x,y]_{\mathfrak a})} =\overline{(0,[fx,u-fy]_{\mathfrak b})}.
\end{equation*}
Hence the quotient identification $(y,u)+G_f\mapsto u-fy$ carries this action to $u\mapsto[fx,u]_{\mathfrak b}$. Together with $(x,fx)\mapsto x$, it therefore identifies the Chevalley--Eilenberg differentials term by term in every degree. Under these identifications, the first normal differential identity becomes exactly \eqref{eq:differential-comparison}.
\end{proof}


\begin{thebibliography}{99}

\bibitem{CrainicSchaetzStruchiner}
M. Crainic, F. Sch\"atz and I. Struchiner,
A survey on stability and rigidity results for Lie algebras,
\emph{Indag. Math. (N.S.)} \textbf{25} (2014), 957--976.

\bibitem{FregierZambon}
Y. Fr\'egier and M. Zambon,
Simultaneous deformations of algebras and morphisms via derived brackets,
\emph{J. Pure Appl. Algebra} \textbf{219} (2015), 5344--5362.

\bibitem{Getzler}
E. Getzler,
Lie theory for nilpotent $L_\infty$-algebras,
\emph{Ann. of Math. (2)} \textbf{170} (2009), 271--301.

\bibitem{Hajjaji}
A. Hajjaji,
Maurer--Cartan characterizations and cohomologies of crossed homomorphisms on Lie triple systems,
\emph{Comm. Algebra} \textbf{52} (2024), 825--844.

\bibitem{ChtiouiHajjajiMabroukMakhlouf}
T. Chtioui, A. Hajjaji, S. Mabrouk and A. Makhlouf,
Cohomologies and deformations of $\mathcal O$-operators on Lie triple systems,
\emph{J. Math. Phys.} \textbf{64} (2023), 081701.

\bibitem{Harris}
B. Harris,
Cohomology of Lie triple systems and Lie algebras with involution,
\emph{Trans. Amer. Math. Soc.} \textbf{98} (1961), 148--162.

\bibitem{Jacobson}
N. Jacobson,
Lie and Jordan triple systems,
\emph{Amer. J. Math.} \textbf{71} (1949), 149--170.

\bibitem{JiangShengSun}
J. Jiang, Y. Sheng and G. Sun,
Stability and rigidity of $3$-Lie algebra morphisms,
\emph{Differential Geom. Appl.} \textbf{101} (2025), 102278.

\bibitem{KuboTaniguchi}
F. Kubo and Y. Taniguchi,
A controlling cohomology of the deformation theory of Lie triple systems,
\emph{J. Algebra} \textbf{278} (2004), 242--250.

\bibitem{Lister}
W. G. Lister,
A structure theory of Lie triple systems,
\emph{Trans. Amer. Math. Soc.} \textbf{72} (1952), 217--242.

\bibitem{MondalSaha}
B. Mondal and R. Saha,
Deformation cohomology of morphisms of Lie--Yamaguti algebras,
\emph{J. Algebra} \textbf{658} (2024), 660--685.

\bibitem{Smirnov}
O. N. Smirnov,
Imbedding of Lie triple systems into Lie algebras,
\emph{J. Algebra} \textbf{341} (2011), 1--12.
\bibitem{Voronov}
Th. Voronov,
Higher derived brackets and homotopy algebras,
\emph{J. Pure Appl. Algebra} \textbf{202} (2005), 133--153.

\bibitem{XiaShengTang}
H. Xia, Y. Sheng and R. Tang,
Cohomology and homotopy of Lie triple systems,
\emph{Comm. Algebra} \textbf{52} (2024), 3622--3642.

\bibitem{Yamaguti}
K. Yamaguti,
On the cohomology space of Lie triple system,
\emph{Kumamoto J. Sci. Ser. A} \textbf{5} (1960), 44--52.

\bibitem{Zhang}
T. Zhang,
Notes on cohomologies of Lie triple systems,
\emph{J. Lie Theory} \textbf{24} (2014), 909--929.

\end{thebibliography}
\end{document}